\documentstyle[12pt,twoside, epsf]{article}

\newcommand{\semi}{\mbox{$\times\!\rule{.2mm}{2mm}$}}
\newcommand{\NN}{\mbox{$I\!\! N$}}	
\newcommand{\ZZ}{\mbox{$Z\!\!\! Z\!$}}	
\newcommand{\QQ}{\mbox{$Q\!\!\!\! I$}}	
\newcommand{\HH}{ {\cal H}_n(q,\infty) }

\newcommand{\BB}{ B_{1,n} }
\newcommand{\WW}{ W_{n,\infty} }    
\newenvironment{pf}{\noindent {\em Proof.}}{\hfill $\Box$ \medskip}

\newtheorem{th}{Theorem}
\newtheorem{lem}{Lemma}

\newtheorem{propn}{Proposition}

\newtheorem{defn}{Definition}

\newtheorem{rem}{Remark}      
\newtheorem{note}{Note}

\def\picill#1by#2(#3)
{\vbox to #2
{\hrule width #1 height 0pt depth 0pt
\vfill\epsffile{#3}}}

\let \ttorg \tt \def \tt{\ttorg \obeyspaces}

\begin{document}
\pagestyle{myheadings}

\markboth{{\sc S. Lambropoulou}}{{\sc Knot theory and $\cal B$-type Hecke algebras}}

\title{Knot theory related to generalized and cyclotomic Hecke algebras of type $\cal B$ }

\author{Sofia Lambropoulou  \\
Mathematisches Institut, Universit\"at G\"ottingen}

 \date{}
\maketitle


\section{Introduction}

After Jones's construction of the classical by now Jones polynomial for knots in $S^3$ using
Ocneanu's Markov trace on the associated Hecke algebras of type $\cal A$, arised questions about
similar constructions on other Hecke algebras as well as in other 3-manifolds. 
\bigbreak
In \cite{L1} is
established that knot isotopy in a 3-manifold may be interpreted in terms of Markov braid
equivalence and, also, that the braids related to the 3-manifold form algebraic structures.
Moreover, the sets of braids related to the solid torus or to the lens spaces $L(p,1)$ 
form groups, which   are in fact the Artin braid groups of type $\cal B$. As a consequence, in 
\cite{L1,L2} appeared the first construction of a Jones-type invariant using
Hecke algebras of type $\cal B$, and this had a natural interpretation
as an isotopy invariant for oriented knots in a solid torus. In a further `horizontal'
development and using a different technique we constructed in \cite{GL} all such solid torus
knot invariants derived from the Hecke algebras of type $\cal B$. Furthermore, in \cite{G} all
Markov traces related to the Hecke algebras of type $\cal D$ were consequently constructed.   
\bigbreak
In this paper we consider  all possible generalizations of the
$\cal B$-type Hecke algebras, namely the cyclotomic and what we call 'generalized', and we construct
Markov traces on each of them, so as to obtain {\em  all possible different levels} of homfly-pt
analogues in the solid torus related to  the  (Hecke) algebras of $\cal B$-type. Our strategy is based
on the one in \cite{L2}, which in turn followed \cite{J}.  So, in this sense, the construction in
\cite{L1,L2} is incorporated here as the most basic level. 

\bigbreak 
In more detail: It is well-understood from Jones's construction of the
homfly-pt (2-variable Jones) polynomial, $P_L$, in \cite{J}, that ${\cal H}_n(q)$, the
Iwahori-Hecke algebra of ${\cal A}_n$-type, is  a quotient of the braid group algebra  $\ZZ
\,  [q^{\pm 1}]B_n$ by factoring out the quadratic relations  
\[ \sigma_i^2=(q-1)\sigma_i+q \]  
\noindent and that these relations  reflect precisely the skein property of $P_L$: 

\[\frac{1}{\sqrt{q}\sqrt{\lambda}}\, P_{L_+} - \sqrt{q}\sqrt{\lambda}\,
 P_{L_-}  = (\sqrt{q}-\frac{1}{\sqrt{q}})\,P_{L_0},  \] 

\noindent where  $L_+$ is a regular projection of an oriented link containing a specified positive
crossing, $L_-$ the same projection with a negative crossing instead, and $L_0$ yet the same
projection with no crossing.    
\bigbreak
We do now analogous considerations  for the solid torus, which we denote by $ST$.  Let us
consider the following Dynkin diagram.

\begin{picture}(300,50)
\put(  0, 20){$({\cal B}_{n})$}
\put( 40, 20){\circle*{5}}
\put( 40, 18){\line(1,0){20}}
\put( 40, 22){\line(1,0){20}}
\put( 60, 20){\circle*{5}}
\put( 60, 20){\line(1,0){30}}
\put( 90, 20){\circle*{5}}
\put( 90, 20){\line(1,0){20}}
\put(120, 20){\circle*{1}}
\put(130, 20){\circle*{1}}
\put(140, 20){\circle*{1}}
\put(150, 20){\line(1,0){20}}
\put(170, 20){\circle*{5}}
\put( 37, 30){$t$}
\put( 57, 30){$\sigma_1$}
\put( 87, 30){$\sigma_2$}
\put(167, 30){$\sigma_{n-1}$}
\put(250, 20){$n \geq 1$}
\end{picture}

\noindent The symbols $t,\sigma_1,\ldots,\sigma_{n-1}$ labelling the nodes correspond to the 
generators of  the  Artin braid group of type ${\cal B}_n$, which we denote by $B_{1,n}$. 
$B_{1,n}$ is defined therefore by the  relations 
\[\begin{array}{rclll}
\sigma_1t\sigma_1t & = & t\sigma_1t\sigma_1  & & \\
t\sigma_i & = & \sigma_it & \mbox{ if } & i>1 \\
\sigma_i\sigma_j & = & \sigma_j\sigma_i & \mbox{ if } & |i-j|>1 \\
\sigma_i\sigma_{i+1}\sigma_i & = & \sigma_{i+1}\sigma_i\sigma_{i+1}  & \mbox{ if } & 1 \leq i
\leq n-2 \end{array}\]
Relations of these  types will  be called {\it  braid relations}.  

$B_{1,n}$ may be seen as the 
subgroup of $B_{n+1}$, the classical braid group on $n+1$ strands, the elements of which keep the
first strand fixed (this is the reason for having chosen the symbol $B_{1,n}$).   This allows
for a geometric interpretation of the elements of $B_{1,n}$ as  mixed braids in $S^3$.  Below
we illustrate the generators $\sigma_i, t$ and the element $t_i'=\sigma_i \ldots
\sigma_1t\sigma_1^{-1} \ldots \sigma_i^{-1}$ in $B_{1,n}$, which plays a crucial role in this
work. 

 $$\vbox{\picill3.7inby1.2in(canon1)  }$$

\noindent Note that the inverses of $\sigma_i, t$ are represented by the same geometric 
pictures, but with the  opposite crossings. 

As shown in \cite{L1,L2}, we can represent oriented knots and links inside $ST$ by
elements of the groups $B_{1,n}$, where the fixed strand represents the complementary solid
torus in $S^3$, and the next $n$ numbered strands represent the 
knot in $ST$. Also, that knot
isotopy in $ST$ can be translated in terms of equivalence classes in 
$\bigcup_{n=1}^{\infty}B_{1,n}$ (Markov theorem), the equivalence being generated by the following
two moves.
 \begin{itemize}
\item[(i)] Conjugation: if $\alpha,\beta \in B_{1,n}$
then $\alpha \sim \beta ^{-1}\alpha\beta$. 
\vspace{-.1in}
\item[(ii)] Markov moves: if $\alpha \in B_{1,n}$ then $\alpha \sim
\alpha{\sigma_n}^{\pm 1}\in B_{1,n+1}$. 
\end{itemize} 

Consider now the classical Iwahori-Hecke algebra of type ${\cal B}_n$, ${\cal H}_n(q,Q)$, as a
quotient of the group algebra  $\ZZ\, [q^{\pm 1}, Q^{\pm 1}] B_{1,n}$  by factoring
out the ideal generated by the  relations $ t^2=(Q-1)t+Q $  and $ g_i^2=(q-1)g_i+q 
$ for all  $i$,  where we denote the image of $\sigma_i$ in  ${\cal H}_n(q,Q)$  by $g_i$. The
idea  in \cite{L1,L2,GL} was to construct invariants  of knots in the solid torus by
constructing trace  functions  $\tau$ on $ \bigcup_{n=1}^{\infty} {\cal H}_n(q,Q)$ which 
support the Markov property:  
 \[ \tau(hg_n)=z\tau(h),\] 
\noindent for $z$  an independent variable in $\ZZ\, [q^{\pm 1}, Q^{\pm 1}]$
and   $h \in {\cal H}_n(q,Q)$. In other words, traces that respect  the above braid
equivalence on  $\bigcup_{n=1}^{\infty}B_{1,n}$. The construction of such traces was only possible
because we were able to find an appropriate inductive basis on $H_{n+1}(q,Q)$, every
element of which  involves the generator
$g_n$ or  the  element $t_n':=g_n \ldots g_1tg_1^{-1} \ldots g_n^{-1}$ at most once (see picture
above for the lifting of $t_i'$ in $B_{1,n}$). In particular, the trace constructed in
\cite{L1,L2} was well-defined inductively by the rules: 
\[ \begin{array}{lll} 
1) & tr(ab)=tr(ba)    & a,b \in {\cal H}_n(q,Q) \\  
2) & tr(1)=1           & \mbox{for all }  {\cal H}_n(q,Q) \\ 
3) & tr(ag_n)=z\, tr(a) & a \in {\cal H}_n(q,Q) \\ 
4) & tr(at_n')=s\, tr(a) & a \in {\cal H}_n(q,Q) 
\end{array} \]

\noindent If we had not used the elements $t_n'$ in the above constructions we would have not been
able to define the trace with only four simple rules. The intrinsic reason for this is that
$B_{1,n}$ splits as a semi-direct product of the classical braid group $B_n$ and of its free 
subgroup $ P_{1,n} $ generated precisely by the elements $t, t_1', \ldots, t_{n-1}'$: 

\[ B_{1,n} = P_{1,n} \semi \, B_n.\]

The Jones-type invariants in $ST$ constructed from the above traces on  $ \bigcup_{n=1}^{\infty}
{\cal H}_n(q,Q)$ satisfy the skein rule related to the quadratic relations  $ g_i^2=(q-1)g_i+q 
$ plus another one reflecting the quadratic relation $ t^2=(Q-1)t+Q $  (cf.
\cite{L1,L2,GL} for an extensive treatment). 

\bigbreak
During the work of S.L. and J. Przytycki on the problem of computing the 3rd skein module of the lens
spaces $L(p,1)$ following the above strategy, it turned out that the skein rule of the homfly-pt
type invariants in \cite{L1,L2,GL} related to $t$ was actually 'artificial', so far that  knot
invariants in $ST$ were concerned, and that for analogous constructions in  $L(p,1)$  it was  needed 
to have constructed first the most generic 2-variable Jones analogue in $ST$, one that would not
satisfy any skein relation involving $t$.  \bigbreak We drop then the  quadratic relation of $t$, and
we consider  the quotient  of the group algebra  $\ZZ\, [q^{\pm 1}] B_{1,n}$  by factoring out only
the relations   $$ g_i^2=(q-1)g_i+q $$
\noindent for all  $i$. This is now a new infinite dimensional algebra, which we denote by ${\cal
H}_n(q, \infty)$ and we shall call it {\em generalized Iwahori-Hecke algebra of type $\cal B$}.
By  $g_i$ above  we denote the image of $\sigma_i$ in  ${\cal H}_n(q, \infty)$, whilst the
symbol  $\infty$ was chosen to indicate that the generator $t$ satisfies no order relation  (since
now  any power  $t^k,$ \ for $ k\in \ZZ \, $ may appear, like in $ B_{1,n}$). For connections of 
these algebras with the affine Hecke algebras of type $\cal A$ see Remark 1. 

\bigbreak
But we would like now to go one step back and, instead of removing from  
${\cal H}_n(q,Q)$ the quadratic relation for $t$, to require that $t$ satisfies a relation
given by a cyclotomic polynomial  of degree $d$:

\[(t-u_1)(t-u_2)\cdots (t-u_d)=0 \]

\noindent Then we  obtain a finite-dimensional algebra known as {\em  cyclotomic Hecke
algebra  of type $\cal B$}, denoted here by ${\cal H}_n(q, d)$. The corresponding {\em 
cyclotomic  Coxeter group of type $\cal B$}, which we 
denote by $W_{n,d}$, is obtained as a quotient of $B_{1, n}$ modulo the relations $g_i^2=1$
and $t^d=1$. ${\cal H}_n(q, d)$ may be seen as a `$d$-deformation' of $W_{n,d}$: In order to
obtain the group algebra we have to substitute the parameters of the cyclotomic polynomial by
the $d$ th roots of unity (and not by $1$ as in the classical case).   These  algebras 
have been introduced  and  studied independently by two groups of mathematicians in
\cite{A2,AK,BM,BkM}.  It follows from the  discussion above that the cyclotomic Hecke algebras
are  also related to the knot theory of the solid torus and, in fact, they make the bridge
between ${\cal H}_n(q,Q)$ and ${\cal H}_n(q, \infty)$.

\bigbreak

Like for the classical Hecke algebras of type $\cal B$, in order to construct linear Markov
traces on  $\bigcup_{n=1}^{\infty} {\cal H}_n(q,\infty)$ or on $\bigcup_{n=1}^{\infty} {\cal
H}_n(q,d)$, we need to find appropriate inductive bases on both types of these algebras. The
inductive bases are derived from known basic sets. This is the aim and the main result of Section 3.
Note that, in the  case of ${\cal H}_n(q,Q)$,  we could easily yield such an inductive basis using 
the results in \cite{DJ}, whilst for ${\cal H}_n(q,d)$ we use the results in \cite{AK}, \cite{BM}.
For ${\cal H}_n(q,\infty)$  we  study its structure in Section 2  and we construct a basis for it
using the structure of the braid group $ B_{1,n}$ and the known bases for ${\cal H}_n(q,d)$. 

\bigbreak
\noindent In Section 4 we construct Markov traces on  
$\bigcup_{n=1}^{\infty} {\cal H}_n(q,\infty)$ and on $\bigcup_{n=1}^{\infty} {\cal H}_n(q,d)$
using the inductive bases of Section 3. Finally in Section 5, we normalize the traces according to
the Markov braid theorem in order to  derive the corresponding knot invariants in $ST$, and we
also give skein interpretations. The invariant related to  ${\cal H}_n(q,\infty)$ is the most
interesting one for us, and in this sense, this work may be seen as the required fundament for
extending such constructions to knots in the lens spaces (see remarks at the end). In the special 
case of  ${\cal H}_n(q,\infty)$ the derived knot invariant reproves the  structure of the 3rd skein
module of the solid torus (cf. \cite{HK,T}). On the other hand, the knot invariants
derived from ${\cal H}_n(q,d)$ are related to submodules of the 3rd skein module of $ST$. It may be
worth noting that introducing and studying  ${\cal H}_n(q,\infty)$ has been independent of the
studies on the cyclotomic analogues.

\bigbreak

Our method shows on one hand that the original strategy of \cite{J} can carry through to so
complicated structures. On the other hand it unifies the construction for all these different
${\cal B}$-type algebras and it highlights the algebraic background underlying these knot
invariants in $ST$. The tedious calculations employed for constructing appropriate bases reflect the
tedious arguments of a more combinatorial approach.  

\bigbreak
 It gives the  author pleasure  to acknowledge her thanks to  V.F.R. Jones for his
valuable comments on this work and to  T.~tom Dieck for discussions and valuable
suggestions. Many thanks are also due to M.~Geck for discussions, useful comments and for pointing
out the literature on the cyclotomic Hecke algebras of type $\cal B$, and especially to 
J.~Przytycki for our discussions on the structure of the generalized Coxeter
groups and Hecke algebras.  Finally, financial support by the SFB 170 in G\"{o}ttingen and the
European Union for parts of this work are gratefully acknowledged. 

\section{Finding a basis for ${\cal H}_n(q, \infty)$}

We start by introducing in more detail ${\cal H}_n(q, \infty)$, ${\cal H}_n(q, d)$ and  their
corresponding Coxeter-type groups $W_{n, \infty}$, $W_{n, d}$. 

\begin{defn} \rm The generalized Iwahori-Hecke algebra of type ${\cal B}_n$ is defined as 

\[ {\cal H}_n(q,\infty) := \ZZ\, [q^{\pm 1}] \, B_{1,n}\, /  <{\sigma_i}^2=(q-1) \, \sigma_i +
q \mbox{ \ for all } i >.\]
\noindent The underlying  generalized Coxeter group of type ${\cal B}_n$ is defined as

\[ \WW := B_{1,n}\, /  <{\sigma_i}^2=1 \mbox{ \ for all } i >.\]
 
\end{defn}

\noindent It follows that if $g_i$ denotes the image of $\sigma_i$ in ${\cal H}_n(q, \infty)$,
then  ${\cal H}_n(q,\infty)$ is defined by the generators $ t, g_1, g_2, \ldots ,
g_{n-1} $ and their relations: 

\[\begin{array}{rclll}
tg_1tg_1 & = & g_1tg_1t  & & \\
tg_i & = & g_it & \mbox{ for } & i>1 \\
g_ig_{i+1}g_i & = & g_{i+1}g_ig_{i+1}  & \mbox{ for } & 1 \leq i \leq n-2 \\
g_i g_j & = & g_j g_i & \mbox{ for } & |i-j|>1 \\
{g_i}^2 & = & (q-1) \, g_i +q & \mbox{ for } & \mbox{all $i$ } 
 \end{array}\]
\bigbreak
\noindent ${\cal H}_n(q,\infty)$ is an associative algebra with $1$. Also, it is easily verified
that, if $S_n$ is the symmetric group, then 

$$ \WW ={\ZZ  }^ n \, \semi \, S_n \ (\mbox{ compare with the structure of} B_{1,n}).$$

\begin{defn} \rm   Let ${\cal R} := \ZZ \, [q^{\pm 1}, u_1^{\pm 1}, \ldots, u_d^{\pm 1}, \ldots]$,
where $q, u_1, \ldots, u_d, \ldots$ are indeterminates. The cyclotomic Iwahori-Hecke algebra of
type ${\cal B}_n$ and of degree $d$ is defined as 

\[ {\cal H}_n(q,d) := {\cal R} \, B_{1,n}\, /  <{\sigma_i}^2=(q-1) \, \sigma_i + 
q \mbox{ \ all } i, \, (t-u_1)(t-u_2)\cdots(t-u_d)=0 >.\]

\noindent The underlying cyclotomic Coxeter group of type $\cal B$ and of degree $d$ is:  

\[ W_{n,d} := B_{1,n}\, /  <{\sigma_i}^2=1 \mbox{ \ for all } i, \  t^d=1, \,   d\in \NN 
>.\] 

\end{defn}
\bigbreak
\noindent  The relation $t^d=1$ is derived by the cyclotomic polynomial by substituting the
$u_i$'s by the $d$'th roots of unity. Also, the Coxeter group of ${\cal B}_n$-type, in our
notation $W_{n,2}$, is the quotient of  $B_{1,n}$ over the relations $t^2 = {\sigma_i}^2 = 1$, \ 
for all $i$.  \medskip

\noindent ${\cal H}_n(q,d) $ is  an associative algebra with $1$, and it is a free module over
${\cal R}$ of rank $d^n\cdot n!$, which is precisely the order of $W_{n,d}$ (cf.
\cite{AK},\cite{BM}). If $d=1$ and $u_1=1$, then ${\cal H}_n(q,1)$ is isomorphic to the
Iwahori-Hecke algebra of type ${\cal A}$ (over  $\ZZ \, [q^{\pm 1}]$). If $d=2$, $u_1=-1$ and
$u_2=Q$, we recover the familiar relation of ${\cal H}_n(q,Q)$, the Iwahori-Hecke algebra of
type ${\cal B}$ (over  $\ZZ \, [q^{\pm 1}, Q^{\pm 1}]$). In ${\cal H}_n(q,d)$ we have 

\[\spadesuit \ \ \ \ \ \ \ \  t^d=a_{d-1}t^{d-1}+\cdots + a_0, \ \ \mbox{ where } \] 

\noindent $ a_{d-1}=u_1+\cdots +u_d, \ a_{d-2}=-(u_1u_2+\cdots +u_{d-1}u_d), \ \ldots, 
 \ a_0=(-1)^d(u_1\ldots u_d)$; from this we can derive easily a relation for $t^{-1}$.

\bigbreak

\noindent $W_{n,d}$ may also be seen as the quotient $\WW \, / <t^d=1>$, \   $d\in \NN$ of
$\WW$, and it is easily verified that

\[ W_{n,d} ={\ZZ_d} ^{n} \, \semi \, S_n\]
 
\noindent Its order is $d^n\cdot n!,$ whilst  $W_{n,2} = {\ZZ_2}^n\, \semi \, S_n$ (compare with
the structure of $B_{1,n}$).

\begin{note} 
\rm W.l.o.g. we extend the ground ring of  ${\cal H}_n(q, \infty)$ to ${\cal R}$. Then 
${\cal H}_n(q,d) $ may also be obtained from  ${\cal H}_n(q, \infty)$ by factoring out the
cyclotomic relation. In this sense ${\cal H}_n(q,d) $ is a `bridge' between  ${\cal
H}_n(q,\infty)$ and ${\cal H}_n(q,Q)$, the classical Hecke algebra. 
\end{note}

\bigbreak
We shall  now find  a basis for ${\cal H}_n(q,\infty)$ as follows: We  find first a canonical
form for the braid group $B_{1,n}$, which yields a basis for $\ZZ \, [q^{\pm 1}] B_{1,n}$. The
 images of these basic  elements in ${\cal H}_n(q,\infty)$ through the canonical map span  ${\cal
H}_n(q,\infty)$. In \cite{AK,BM} bases for ${\cal H}_n(q,d) $ have been constructed. We then treat
the spanning set and using these bases we obtain a basis for ${\cal H}_n(q,\infty)$. This 
approach 
shows clearly the relation among the structures of $B_{1,n}$, ${\cal H}_n(q,\infty)$,
 ${\cal H}_n(q,d)$ and  $W_{n,\infty}$, $W_{n,d}$.

\bigbreak
\noindent In order to proceed we need to recall the notion of the pure 
braid group and Artin's canonical form for pure braids:  The classical pure braid group, $P_n$,
consists of all elements in $B_n$ that induce the identity permutation in $S_n$; $P_n \lhd B_n$
and $P_n$ is generated by the elements

\vspace{.1in}
$ \begin{array}{ll}
A_{rs} & =  {\sigma_r}^{-1}{\sigma_{r+1}}^{-1}\ldots
{\sigma_{s-2}}^{-1}{\sigma_{s-1}}^2 \sigma_{s-2}\ldots\sigma_{r+1}\sigma_r \\ [.1in] 
  & =  \sigma_{s-1}\sigma_{s-2}\ldots \sigma_{r+1}{\sigma_r}^2
{\sigma_{r+1}}^{-1}\ldots {\sigma_{s-2}}^{-1}{\sigma_{s-1}}^{-1}, \ \  1\leq r<s\leq n. 
\end{array} $

\bigbreak
\noindent  Artin's canonical form says that every element, $A$, of $P_n$ can be written uniquely
in the form:  
\[ A=U_1U_2 \cdots U_{n-1} \] 

\noindent where each $U_i$ is a uniquely determined product of powers of the $A_{ij}$ using only
those with $i<j$. Geometrically, this means that any pure braid can be `combed' i.e. can be written
 canonically as: the pure braiding of the first string with the rest, then keep the first string
fixed and uncrossed and have the pure braiding of the second string and so on (cf. [J.S.~Birman, 
Braids, Links and Mapping Class Groups,  {\it Ann. of Math. Stud.}  82, Princeton University Press,
Princeton 1974] for a complete treatment). 

\bigbreak

We find now a canonical form for $B_{1,n}$.  An element $w$ of $ B_{1,n}$ induces a
permutation $\sigma \in S_n$ of the $n$ numbered strands. We add at the bottom of the braid a
standard braid in $B_n$ corresponding to $\sigma ^{-1}$, and then we add its inverse $\sigma $. 
 Now, $w\sigma^{-1}$ is a pure braid  on $n+1$ stands (including the first fixed one),
and we apply to it Artin's canonical form. This separates the braiding of the fixed strand from
the rest:

$$\vbox{\picill5inby1.6in(canon4)  }$$

The above is in fact  the proof of the decomposition of $B_{1,n}$ as a semidirect product:

\begin{propn} $B_{1,n} = P_{1,n} \semi \, B_n$. 
\end{propn}

\noindent From the uniqueness of Artin's canonical form, it
follows that any $w \in B_{1,n}$ can be expressed uniquely as a product $ v\cdot \sigma$  \
(`vector-permutation'), where $v$ is an element of the free group $ P_{1,n}:$

\[ v={t'_{i_1}}^{k_1}{t'_{i_2}}^{k_2}\ldots {t'_{i_r}}^{k_r}, \  k_1, \ldots, k_r \in \ZZ,  
\mbox{ \ where \ } {t'_i}^k := \sigma_i \ldots \sigma_1t^k\sigma_1^{-1} \ldots \sigma_i^{-1}, \]

\noindent and $\sigma \in B_n$ is written in the induced by $P_n$
canonical form. Thus the set $\{ v\cdot \sigma \}$ forms a basis for the algebra $\ZZ \, [q^{\pm
1}] \, B_{1,n}$, and, therefore, it spans the quotient ${\cal H}_n(q,\infty)$. On the level of 
${\cal H}_n(q,\infty)$ we can already improve this spanning set, since on this level  $\sigma$ is
a word in  ${\cal H}_n(q)$, the Iwahori-Hecke algebra of ${\cal A}_{n-1}$-type. So, $\sigma$ can
be written in terms of the standard basis of ${\cal H}_n(q)$ (cf. \cite{J}):  
 
\[ \{ (g_{i_1}g_{i_1-1}\ldots g_{i_1-r_1})(g_{i_2}g_{i_2-1}\ldots g_{i_2-r_2})\ldots
(g_{i_p}g_{i_p-1}\ldots g_{i_p-r_p}) \}, \] 
\[ \mbox{ for \ } 1\leq i_1<\ldots <i_p\leq n-1 \ \mbox{ and \ } r_j\in\{0, 1, \ldots, i_j-1\}.\]

\noindent Therefore we showed

\begin{propn} The set 

\[ {\Sigma}_1 = \{ {t'_{j_1}}^{k_1}{t'_{j_2}}^{k_2}\ldots {t'_{j_r}}^{k_r} \cdot \sigma \} , \] 

\noindent where $ t'_{0}:=t, \ {t'_i}^k := g_i \ldots g_1t^kg_1^{-1} \ldots g_i^{-1}, \    
j_1,\ldots, j_r \in \{ 0, 1, \ldots, n-1 \}$, \ $ k_1, \ldots, k_r \in \ZZ $  \ and \ $\sigma$  a
basic element of ${\cal H}_n(q)$,   spans \ ${\cal H}_n(q,\infty).$ 
\end{propn}
\noindent Notice that the indices of the `vector' part are not ordered. Also, that the above
canonical form for $\BB$ yields immediately the following canonical form  $\{ v\cdot \sigma \}$
for  $\WW$:   

\[ \{ v\cdot \sigma \} = \{ {t_{j_1}}^{k_1}{t_{j_2}}^{k_2}\ldots {t_{j_r}}^{k_r} \cdot
\sigma \} , \] 

\noindent where $ t_{0}:=t, \ {t_i}^k := s_i \ldots s_1t^ks_1 \ldots s_i,$ \ for \ 
$ 0\leq j_1<\ldots <j_r\leq n-1$, \ $ k_1, \ldots, k_r \in \ZZ $  \ and  \ $\sigma \in S_n$ is an
element of the canonical form of $S_n$ \ (where $s_i$ denotes the image of $\sigma_i$ in $\WW$).
Thus,  this set also forms a  basis for the group algebra $\ZZ \, [q^{\pm 1}] \, \WW$. 
\bigbreak

\noindent   Notice
here that the indices of the `vector' part are ordered. This suggests that it  may be possible to
order the indices $j_1,\ldots, j_r$ of the words  ${t'_{j_1}}^{k_1}{t'_{j_2}}^{k_2}\ldots
{t'_{j_r}}^{k_r}$ in ${\Sigma}_1$, so as to be left with a canonical basis for $\HH$.   To
achieve this straight from ${\Sigma}_1$ is  very difficult, because it is hard to get hold of
an induction step, even though there are relations among the ${t'_i}^{k_i}$'s. Instead, we change
the ${t'_i}^k$'s to the elements ${t_i}^k$, where $t_{0}:=t,$ \ and \ $t_i := g_i \ldots g_1tg_1
\ldots g_i $. These elements  commute in $\HH$.

\bigbreak
The following relations hold in $\HH$ and in ${\cal H}_n(q,d)$  and 
will be used repeatedly in the sequel.

\begin{lem} For $\epsilon \in \{ \pm 1 \}$ the following hold:
\begin{itemize}
\item[(i)]  $ {g_i}^\epsilon = q^\epsilon \, {g_i}^{- \epsilon} +(q^\epsilon -1),$
   
${g_i}^{2 \epsilon}=(q^\epsilon -1) \, {g_i}^\epsilon + q^\epsilon , \ \mbox{ for } q\neq 0. $

\item[(ii)] $ {g_i}^\epsilon ({g_k}^{\pm 
1}g_{k-1}^{\pm 1}\ldots {g_j}^{\pm 1}) = ({g_k}^{\pm 1}g_{k-1}^{\pm 1}\ldots {g_j}^{\pm
1}){g_{i+1}}^\epsilon, \ \mbox{ for } k>i\geq j,$ 

${g_i}^\epsilon ({g_j}^{\pm
1}g_{j+1}^{\pm 1}\ldots {g_k}^{\pm 1}) = ({g_j}^{\pm 1}g_{j+1}^{\pm 1}\ldots {g_k}^{\pm
1}){g_{i-1}}^\epsilon, \ \mbox{ for } k\geq i> j, $

where the sign of the ${\pm 1}$ superscript is the same for all generators.

\item[(iii)] $ g_ig_{i-1}\ldots g_{j+1}{g_j} g_{j+1}\ldots g_i =
g_jg_{j+1}\ldots g_{i-1}{g_i} g_{i-1}\ldots g_{j+1}{g_j},$ 

$ {g_i}^{- 1}g_{i-1}^{- 1}\ldots g_{j+1}^{- 1}{g_j}^\epsilon g_{j+1}\ldots g_i =
g_jg_{j+1}\ldots g_{i-1}{g_i}^\epsilon g_{i-1}^{- 1}\ldots g_{j+1}^{- 1}{g_j}^{- 1}.$  

\item[(iv)]  $ {g_i}^\epsilon \ldots {g_{n-1}}^\epsilon {g_n}^\epsilon {g_n}^\epsilon
{g_{n-1}}^\epsilon\ldots {g_i}^\epsilon= $

$ (q^\epsilon -1) \, \sum_{r=0}^{n-i} \, q^{\epsilon r} \, ({g_i}^\epsilon 
\ldots {g_{n-r-1}}^\epsilon {g_{n-r}}^\epsilon {g_{n-r-1}}^\epsilon  \ldots
{g_i}^\epsilon)+q^{\epsilon (n-i+1)}=$

$ \sum_{r=0}^{n-i+1} \, (q^\epsilon -1)^{\epsilon_r} q^{\epsilon r} \, ({g_i}^\epsilon
 \ldots {g_{n-r-1}}^\epsilon {g_{n-r}}^\epsilon {g_{n-r-1}}^\epsilon\ldots {g_i}^\epsilon),$ 

$ \mbox{ where } \ \epsilon_r=1 \ \mbox{ if } \ r\leq n-i \ \mbox{ and } \ \epsilon_{n-i+1}=0.$

Similarly,

${g_i}^\epsilon \ldots {g_2}^\epsilon {g_1}^\epsilon {g_1}^\epsilon 
{g_2}^\epsilon \ldots {g_i}^\epsilon =$ 

$(q^\epsilon-1) \,\sum_{r=0}^{i-1} \, q^{\epsilon r} \, ({g_i}^\epsilon \ldots {g_{r+2}}^\epsilon
{g_{r+1}}^\epsilon {g_{r+2}}^\epsilon \ldots {g_i}^\epsilon)+q^{\epsilon i}=$

$\sum_{r=0}^{i} \, (q^\epsilon -1)^{\epsilon_r} q^{\epsilon r} \, ({g_i}^\epsilon
 \ldots{g_{r+2}}^\epsilon {g_{r+1}}^\epsilon {g_{r+2}}^\epsilon \ldots {g_i}^\epsilon)$, 

$ \mbox{ where } \ \epsilon_r=1 \ \mbox{ if } \ r\leq i-1 \ \mbox{ and } \ \epsilon_i=0$.

\item[(v)]  $t^{\lambda} g_1tg_1 = g_1tg_1 t^{\lambda} \ \mbox{ for } {\lambda}\in \ZZ,$

$g_i{t_k}^\epsilon = {t_k}^\epsilon g_i \ \mbox{ for } \ k>i, \, k<i-1,$

$g_it_i = q \, t_{i-1}g_i+(q-1) \, t_i,$ 

$g_it_{i-1} = q^{-1} \, t_ig_i+(q^{-1}-1) \, t_i,$ 

$g_i{t_{i-1}}^{-1} = q \, {t_i}^{-1}g_i+(q-1) \, {t_{i-1}}^{-1},$
 
$g_i{t_i}^{-1} = q^{-1} \, {t_{i-1}}^{-1}g_i+(q^{-1}-1) \, {t_{i-1}}^{-1}$.

\item[(vi)]  $g_i{t'_k}^\epsilon = {t'_k}^\epsilon g_i \ \mbox{ for } \ k>i, \ k<i-1,$ 

$g_i{t'_i}^\epsilon = {t'_{i-1}}^\epsilon g_i + (q-1) \, {t'_i}^\epsilon + 
(1-q) \, {t'_{i-1}}^\epsilon,$

$g_i{t'_{i-1}}^\epsilon = {t'_i}^\epsilon g_i.$ 
 
\item[(vii)]  ${t_i}^k{t_j}^\lambda = {t_j}^\lambda {t_i}^k \ \mbox{ for } \ i\neq j \ \mbox{ and
} \ k, \lambda\in \ZZ.$

\item[(viii)]  ${t_i'}^k = g_i\ldots g_1 t^k {g_1}^{-1}\ldots {g_i}^{-1} \ \mbox{ for } \ k \in
\ZZ.$ 

Therefore we have in ${\cal H}_n(q,d)$:

$(t'_i - u_1)(t'_i - u_2)\ldots (t'_i - u_d)=0,$  which implies  
${t'_i}^d = a_{d-1}{t'_i}^{d-1}+\cdots +a_0,$

and where the $a_i$'s are given in relation $(\spadesuit)$ in Section 2.

\end{itemize}
\end{lem} 

\begin{pf} 
We point out first that in the rest of the paper and in order to
facilitate the reader  we underline in the proofs the expressions  which are crucial for the
next step. We also use the symbol `$\sum$' instead of the phrase `linear combination of
words of the type'.

Except for (iv), all relations are easy consequenses of the defining relations of
$\HH$ respectively ${\cal H}_n(q,d)$.  Relation (vii) can be also checked using braid diagrams.
We prove (iv) by induction on the length $l=n-i+1$ of the word $g_ng_{n-1}\ldots g_i$. For $l=1$
we have ${g_n}^2=(q-1)g_n+q^1$. Assume now (iv) holds up to $l=n-i$. Then for $l=n-i+1$ we have 

\vspace{.1in}
\noindent $g_i\underline{g_{i+1}\ldots g_ng_n\ldots g_{i+1}}g_i\stackrel{induction \, step}{=}$

\vspace{.1in}
$ \underline{g_i} \, [(q-1) \, \sum_{r=0}^{n-(i+1)} \, q^r \, (g_{i+1}\ldots g_{n-r-1}
g_{n-r} g_{n-r-1}\ldots g_{i+1})+q^{n-i} ] \, \underline{g_i} = $ 

\vspace{.1in}
$ (q-1) \, \sum_{r=0}^{n-(i+1)} \, q^r \, (g_i\ldots g_{n-r-1} g_{n-r} g_{n-r-1}\ldots g_i)
+q^{n-i}\underline{{g_i}^2}= $

\vspace{.1in}
$ (q-1) \, \sum_{r=0}^{n-(i+1)} \, q^r \, (g_i\ldots g_{n-r-1} 
g_{n-r}g_{n-r-1}\ldots g_i) +(q-1)q^{n-i}g_i+q^{n-i+1}= $ 

\vspace{.1in}
$ (q-1) \, \sum_{r=0}^{n-i} \, q^r \, (g_i\ldots g_{n-r-1} g_{n-r} g_{n-r-1}\ldots
g_i)+q^{n-i+1}.$

\noindent Furthermore note that in the Relations (v) and (vi) a $t_i$ or a $t'_i$ will {\it  not}
change to a  ${t_i}^{-1}$ or a ${t'_i}^{-1}$ respectively and, therefore, these relations
preserve the total sum of the exponents of the  $t_i$'s and the $t'_i$'s in a word. Note also
that for $j=i-1$  the relations (iii) boil down to the usual braid relation and its variations
with inverses. 
\end{pf}   

\begin{th}
 In $\HH$ the set 
\[ {\Sigma}_2= \{ {t_{i_1}}^{k_1}{t_{i_2}}^{k_2}\ldots {t_{i_r}}^{k_r}\cdot \sigma \}\] 

\noindent for  $0\leq i_1<\ldots <i_r\leq n-1$, $k_1, \ldots, k_r \in \ZZ$ \ and 
$\sigma$ a basic element in ${\cal H}_n(q)$, forms a basis for ${\cal H}_n(q,\infty).$ 
\end{th}
\noindent  Notice that in ${\Sigma}_2$ the indices of the `vector' part are ordered.  

\vspace{.1in}
\begin{pf} 
To show that 
${\Sigma}_2$ spans $\HH$ it suffices, by Proposition 2, to show that an element of 
${\Sigma}_1$ can be written as a linear combination of elements in ${\Sigma}_2$. Indeed, let 

\[w = {t'_{j_1}}^{k_1}{t'_{j_2}}^{k_2}\ldots {t'_{j_m}}^{k_m}\cdot \sigma \in {\Sigma}_1.\]

\noindent We do the proof by induction on $$\rho = |k_1|+|k_2|+\cdots + |k_m|,$$  
the absolute number of 
$t$'s in $w$. For $\rho = 1$  either $w = t'_i\cdot \sigma \, $ or $ \, w = {t'_i}^{-1}\cdot
\sigma:$

\vspace{.1in}
\noindent $t'_i\cdot \sigma =g_i \ldots g_1t{g_1}^{-1}\ldots {g_i}^{-1}\cdot \sigma=$ 

\vspace{.1in}
$ \underline{g_i\ldots g_1t(g_1\ldots g_i}{g_i}^{-1}\ldots {g_1}^{-1}){g_1}^{-1}\ldots
{g_i}^{-1}\cdot \sigma = t_i\cdot \sigma_1,$ 

\vspace{.1in}
where $\sigma_1 = {g_i}^{-1}\ldots {g_1}^{-1}{g_1}^{-1}\ldots
{g_i}^{-1}\cdot \sigma \in {\cal H}_n(q)$, a linear combination of 

basic elements   of ${\cal H}_n(q)$.

\bigbreak
\noindent ${t'_i}^{-1}\cdot\sigma = g_i \ldots g_1t^{-1}{g_1}^{-1}\ldots {g_i}^{-1}\cdot
\sigma=$ 
    
\vspace{.1in} 
$\underline{g_i\ldots g_1(g_1\ldots g_i}{g_i}^{-1}\ldots
{g_1}^{-1})t^{-1}{g_1}^{-1}\ldots {g_i}^{-1}\cdot \sigma\stackrel{Lemma 1,(iv)}{=}$

\vspace{.1in}
$(q-1) \, \sum_{r=0}^{i-1} \, q^r \, (g_i\ldots g_{r+2}\underline{g_{r+1} g_{r+2}\ldots g_i)
{g_i}^{-1}\ldots {g_{r+1}}^{-1}}{g_r}^{-1}
\ldots {g_1}^{-1}t^{-1}{g_1}^{-1}\ldots $

\vspace{.1in}
${g_i}^{-1}\cdot \sigma + q^i \, t_i^{-1}\cdot \sigma =$

\vspace{.1in}
$(q-1) \, \sum_{r=0}^{i-1} \, q^r \, (g_i\ldots g_{r+2})
\underline{{g_r}^{-1}\ldots {g_1}^{-1}t^{-1}{g_1}^{-1}\ldots{g_r}^{-1}}\ldots {g_i}^{-1}\cdot
\sigma + q^i \, t_i^{-1}\cdot \sigma =$

\vspace{.1in}
$(q-1) \, \sum_{r=0}^{i-1} \, q^r \, {t_r}^{-1}(g_i\ldots g_{r+2}{g_{r+1}}^{-1}\ldots 
{g_i}^{-1}\cdot \sigma) + q^i \, t_i^{-1}\cdot \sigma =$

\vspace{.1in}
$(q-1) \, \sum_{r=0}^{i-1} \, q^r \, {t_r}^{-1}\cdot \sigma_r + q^i \, t_i^{-1}\cdot \sigma,$
 \ where $\sigma_r = g_i\ldots g_{r+2}{g_{r+1}}^{-1}\ldots {g_i}^{-1}\cdot \sigma \in$ 

\vspace{.1in} ${\cal H}_n(q).$

\vspace{.1in}
\noindent  Suppose now the assumption holds for up to $\rho -1$ \ $t$'s in $w$. Then, the
induction step holds in particular for all such words with $\sigma =1$. So, for
$|k_1|+|k_2|+\cdots + |k_m|=\rho $ we have: 

\vspace{.1in}
\noindent ${t'_{j_1}}^{k_1}\ldots {t'_{j_m}}^{k_m}\cdot \sigma = \left\{ 
\begin{array}{ll}
\underline{{t'_{j_1}}^{k_1}\ldots {t'_{j_m}}^{k_m-1}}t'_{j_m}\cdot \sigma, & \mbox{ if } k_m>0
\\ [.1in] 
\underline{{t'_{j_1}}^{k_1}\ldots {t'_{j_m}}^{k_m+1}}{t'_{j_m}}^{-1}\cdot \sigma, 
& \mbox{ if } k_m<0 \end{array} \right. $

\bigbreak
$\stackrel{by \, induction}{=} \left\{  
\begin{array}{ll}
{\Sigma \, t_{i_1}}^{\lambda_1}\ldots {t_{i_n}}^{\lambda_n}\cdot\sigma_1\cdot
\underline{t'_{j_m}} \cdot\sigma, & \mbox{ for some } \sigma_1\in {\cal H}_n(q), \\ [.1in]  
{\Sigma \, t_{\mu_1}}^{\nu_1}\ldots {t_{\mu_n}}^{\nu_n}\cdot\sigma_2 \cdot 
\underline{{t'_{j_m}}^{-1}}\cdot \sigma, & \mbox{ for some } \sigma_2\in {\cal H}_n(q),
\end{array} \right. $

\bigbreak
$\left\{  
\begin{array}{ll}
1\leq i_1<\ldots <i_n\leq n-1, & |\lambda_1|+\cdots + |\lambda _n|=\rho -1  \\ [.1in]
1\leq \mu_1<\ldots <\mu_n\leq n-1, & |\nu_1|+\cdots + |\nu _n|=\rho -1
\end{array} \right. $

\bigbreak
$ = \left\{ 
\begin{array}{l}
{\Sigma \, t_{i_1}}^{\lambda_1}\ldots {t_{i_n}}^{\lambda_n}\cdot\underline{\sigma_1\cdot t_{j_m}}
({g_{j_m}}^{-1}\ldots {g_1}^{-1}{g_1}^{-1}\ldots {g_{j_m}}^{-1})\cdot \sigma \\
[.1in] 
{\Sigma \, t_{\mu_1}}^{\nu_1}\ldots {t_{\mu_n}}^{\nu_n}\cdot\underline{\sigma_2\cdot 
(g_{j_m}\ldots g_1g_1\ldots g_{j_m}){t_{j_m}}^{-1}}\cdot \sigma 
 \end{array}\right.. $

\vspace{.1in}
We apply Lemma 1,(v) on the underlying expressions in order to shift $t_{j_m}$ 

and  ${t_{j_m}}^{-1}$ to the left and we obtain sums of the words:

\bigbreak
$\left\{ 
\begin{array}{ll}
{\Sigma \, t_{i_1}}^{\lambda_1}\ldots {t_{i_n}}^{\lambda_n}\cdot t_{e_1}\cdot \sigma'_1,
 & \sigma'_1 \in {\cal H}_n(q), \ e_1\in \{0, 1, \ldots, n-1\} \\ [.1in] 
{\Sigma \, t_{\mu_1}}^{\nu_1}\ldots {t_{\mu_n}}^{\nu_n}\cdot {t_{e_2}}^{-1}\cdot \sigma'_2,
 & \sigma'_2 \in {\cal H}_n(q), \ e_2\in \{0, 1, \ldots, n-1\}
\end{array}\right. $ 


\bigbreak
$ \stackrel{Lemma 1,(vii)}{=} \left\{ 
\begin{array}{ll}
{\Sigma \, t_{i_1}}^{\lambda_1}\ldots {t_{i_r}}^{\lambda_r} \cdot t_{e_1} \cdot 
{t_{i_{r+1}}}^{\lambda_{r+1}}\ldots {t_{i_n}}^{\lambda_n} \cdot \sigma_i \sigma', 
i_r<e_1<i_{r+1}. \\ [.1in]  
{\Sigma \, t_{\mu_1}}^{\nu_1}\ldots {t_{\mu_k}}^{\nu_k} \cdot
{t_{e_2}}^{-1} \cdot  {t_{\mu_{k+1}}}^{\nu_1}\ldots {t_{\mu_n}}^{\nu_n} \cdot \sigma'_i \sigma, 
\mu_k<e_2<\mu_{k+1}. 
 \end{array}\right. $
\bigbreak
\noindent I.e. in either case we obtained a linear combination of elements of $\Sigma_2.$
\bigbreak
We next show  linear independency of the elements of ${\Sigma}_2$:

\noindent Let  $\sum_{i=1}^{m}
\lambda_i w_i = 0$ for $w_1, w_2, \ldots, w_m \in {\Sigma}_2$. We assume first that the  
exponents of the $t_j$'s in the words $w_i$ are all positive for all $i$, and we choose  $d>k \in
\NN$, where $k$ is the maximum of the exponents of the $t_j$'s in $\sum_{i=1}^{m}\lambda_i w_i$.
Then, the canonical   epimorphism of $\HH$ onto ${\cal H}_n(q,d)$  applied on the equation 
$\sum_{i=1}^{m} \lambda_i w_i = 0$ in  $\HH$ yields the equation  $\sum_{i=1}^{m} \lambda_i w_i =
0$ in ${\cal H}_n(q,d)$. As shown in \cite{AK}, Proposition 3.4 and Theorem 3.10, the
elements of ${\Sigma}_2$ with $0<k_1, \ldots, k_r\leq d-1$ form a basis for ${\cal H}_n(q,d), \,
d\in \NN$. (In \cite{AK} \, $d$ is denoted by $r$, \, ${\cal H}_n(q,d)$ is denoted by  ${\cal
H}_{n,r}$ and $\sigma$ is denoted by $a_w$.)  This implies  $\lambda_i=0, \ i=1, \ldots, m$. 


\vspace{.1in}
Assume finally that  some $w_i$'s contain $t_j$'s with negative exponents. The idea is to resolve
the negative exponents and then refer to the previous case. One way  is to proceed as above, and
after we have projected $\sum_{i=1}^{m}\lambda_i w_i = 0$ on ${\cal H}_n(q,d)$, to resolve the
$t_j$'s with negative exponents using the algebra relations; finally, to conclude  
$\lambda_i=0, \ i=1, \ldots, m$, using induction and arguments from linear algebra.
But we would rather give a more elegant argument, that was suggested by T. tom Dieck. 

\vspace{.1in}
\noindent Namely, let $P$ be the product of all $t_j^k, \, k\in \NN$ for all $j, k$ 
 such that $t_j^{-k}$ is in some $w_i$. Since $P$ is an invertible element of $\HH$,
we have   $\sum_{i=1}^{m}\lambda_i w_i = 0 \Leftrightarrow P \cdot \sum_{i=1}^{m}\lambda_i w_i =
0.$ The last equation is eqivalent to $\sum_{i=1}^{m}  \lambda_i P  w_i = 0$, where the elements
$P  w_i$ are pairwise different and the exponents of the  $t_j$'s  contained in each  $P  w_i$ 
are positive  for all $i$. We then refer to the previous case, and the proof of Theorem 1 is now 
concluded. 

\vspace{.1in}
Thus ${\Sigma}_2$ is a basis of $\HH$, and therefore $\HH$ is a free module. 
\end{pf}   

\begin{rem} \rm 
In \cite{D}, (8.23) tom Dieck establishes an isomorphism between $\HH$ and the twisted
tensor product of the Hecke algebra of the Coxeter group of the affine type 
$\tilde{{\cal A}}_{n-1}$. One can also use the extended affine Hecke algebra of type  
$\tilde{{\cal A}}_{n-1}$ and study quotient maps onto ${\cal H}_n(q,d)$ as defined in     
\cite{A2}, Section 2.1. The same map also works for $\HH$ and it is in fact an isomorphism.

\end{rem}

\section{Inductive bases for ${\cal H}_n(q,\infty)$ and ${\cal H}_n(q,d)$} 

The basis of ${\cal H}_n(q,\infty)$ constructed in the previous section as well as the
corresponding one for ${\cal H}_n(q,d)$ yields an inductive basis for ${\cal H}_n(q,\infty)$
respectively  ${\cal H}_n(q,d)$, which gives rise to another two inductive bases, the last one
being the appropriate for constructing Markov traces on these algebras. Here we give these three
inductive bases and we conclude this section by giving another basic set for  ${\cal
H}_n(q,\infty)$ respectively
 ${\cal H}_n(q,d)$, which is analogous to the set $\Sigma_2$, but using $t'_i$'s instead of 
$t_i$'s. 
\bigbreak 
From now on we shall denote by ${\cal H}_n$ both ${\cal H}_n(q,\infty)$  and ${\cal
H}_n(q,d)$  and  by $W_n$ both $W_{n, \infty}$ and $W_{n, d}$. Also, whenever we refer to 
$k\in \ZZ$ respectively $k\in \ZZ_d$ we shall assume $k\neq 0$. 
We now find the first inductive basis for ${\cal H}_{n+1}$. This on the group level is
an inductive canonical form, and it provides a set of right coset representatives of $W_n$ into 
 $W_{n+1}$, which is completely analogous to \cite{DJ}, p. 456 for ${\cal B}$-type
Coxeter groups.

\begin{lem} For $k\in \ZZ$ the following hold in ${\cal H}_{n+1}(q,\infty)$ 
respectively ${\cal H}_{n+1}(q,d)$:  

\vspace{.1in}
\noindent \ (i) ${t_n}^kg_n \: = \: (q-1) \, \sum_{j=0}^{k-1} \, q^j \, {t_{n-1}}^j{t_n}^{k-j} +
q^k \, g_n {t_{n-1}}^k, \ \ \mbox{ if } k\in \NN \mbox{ and }$

\vspace{.1in}
 $ \: {t_n}^kg_n \: = \: (1-q) \, \sum_{j=0}^{k-1} \, q^j \, {t_{n-1}}^j{t_n}^{k-j} +
q^k \, g_n {t_{n-1}}^k, \ \ \mbox{ if } k\in \ZZ \, - \NN.$ 

\vspace{.1in}
\noindent (ii) ${t_n}^kg_ng_{n-1} \ldots g_i \: = $

\vspace{.1in}
$(q-1) \, \sum_{j=0}^{k-1} \, q^j \,({t_{n-1}}^j g_{n-1}g_{n-2} \ldots g_i){t_n}^{k-j}+$ 

\vspace{.1in}
$ (q-1)q^k \, \sum_{j=0}^{k-1} \, q^j \, ({t_{n-2}}^j g_{n-2}g_{n-3} \ldots g_i)g_n
{t_{n-1}}^{k-j}+$

\vspace{.1in}
$ (q-1)q^{2k} \, \sum_{j=0}^{k-1} \, q^j \, ({t_{n-3}}^j g_{n-3} \ldots
g_i)g_ng_{n-1}{t_{n-2}}^{k-j}$ 

\vspace{.1in}
$+ \cdots +$

\vspace{.1in}
$ (q-1)q^{(n-i)k} \,\sum_{j=0}^{k-1} \, q^j \, ({t_{i-1}}^j) g_ng_{n-1}\ldots g_{i+1}{t_i}^{k-j}$

\vspace{.1in}
$+ q^{(n-i+1)k} \, g_ng_{n-1}\ldots g_i{t_{i-1}}^k, \ \ \mbox{ if } k\in \NN,$

\vspace{.1in}
\noindent whilst for $k\in \ZZ \, - \NN$ we have an analogous formula, only $(q-1)$ is
replaced by $(1-q), \ q^k=q^{-|k|}\mbox{ and } |k-j|+|j|=|k|.$
\end{lem}

\begin{pf} 
We prove (i) for the case $ k>0$ by induction on $k$. (For $k<0$ completely analogous.) For $k=1$
we have  $t_ng_n = (q-1) t_n + q^1 g_n t_{n-1}.$ Suppose the assumption holds for $k-1$. Then for
$k$ we have: 

\vspace{.1in}
\noindent ${t_n}^kg_n = t_n\underline{{t_n}^{k-1}g_n} \stackrel{by \, induction}{=}$

\vspace{.1in}
$\underline{t_n} \, [(q-1) \, \sum_{j=0}^{k-2} \, q^j \, {t_{n-1}}^j{t_n}^{k-1-j} + q^{k-1} \,
g_n {t_{n-1}}^{k-1}] \stackrel{Lemma \, 1,(vii)}{=}$ 

\vspace{.1in}
$(q-1) \, \sum_{j=0}^{k-2} \, q^j \, {t_{n-1}}^j{t_n}^{k-j} + q^{k-1} \, \underline{t_ng_n}
{t_{n-1}}^{k-1} \stackrel{Lemma \, 1,(v)}{=}$ 

\vspace{.1in}
$(q-1) \, \sum_{j=0}^{k-2} \, q^j \, {t_{n-1}}^j{t_n}^{k-j} +
q^{k-1}(q-1) \, \underline{t_n{t_{n-1}}^{k-1}} + q^k \, g_n{t_{n-1}}^k \stackrel{Lemma \,
1,(vii)}{=}$

\vspace{.1in}
$(q-1) \, \sum_{j=0}^{k-1} \, q^j \, {t_{n-1}}^j{t_n}^{k-j} + q^k \, g_n{t_{n-1}}^k .$

\vspace{.1in}
\noindent We prove (ii)  for the case $k>0$ by decreasing induction on $i$. (For $k<0$
completely analogous.)  For $i=n$ we have (i). Assume it holds for $i+1<n \ (\Leftrightarrow
i\leq n-2 \Leftrightarrow n-i\geq 2)$. Then for $i$ we have: 

\vspace{.1in}
\noindent $\underline{{t_n}^kg_n\ldots g_{i+1}} g_i \stackrel{by \, induction}{=}$
\bigbreak
$[(q-1) \, \sum_{j=0}^{k-1} \, q^j \, ({t_{n-1}}^j g_{n-1}g_{n-2} \ldots g_{i+1})
\underline{{t_n}^{k-j}]g_i} + \cdots +$

\vspace{.1in}
$[(q-1)q^{(n-(i+1))k} \, \sum_{j=0}^{k-1} \, q^j \, ({t_i}^j) g_ng_{n-1} \ldots
g_{i+2}\underline{{t_{i+1}}^{k-j}]g_i} + $

\vspace{.1in}
$[q^{(n-i)k} \, g_ng_{n-1} \ldots g_{i+1}\underline{{t_i}^k]g_i} \stackrel{\rm Lemma 1,(v)
\&  Lemma 2,(i)}{ =} $

\vspace{.1in}
$(q-1) \, \sum_{j=0}^{k-1} \, q^j \, ({t_{n-1}}^j g_{n-1}\ldots g_{i+1}g_i){t_n}^{k-j} + \cdots +$

\vspace{.1in}
$(q-1)q^{(n-(i+1))k} \, \sum_{j=0}^{k-1} \, q^j \, ({t_i}^jg_i) g_ng_{n-1} \ldots
g_{i+2}{t_{i+1}}^{k-j} + $

\vspace{.1in}
$q^{(n-i)k}(q-1) \, \sum_{j=0}^{k-1} \, q^j \, g_ng_{n-1} \ldots
g_{i+1}\underline{{t_{i-1}}^j}{t_i}^{k-j}+$

\vspace{.1in}
$q^{(n-i)k}q^k \, g_ng_{n-1} \ldots g_i{t_{i-1}}^k =$

\vspace{.1in}
$(q-1) \, \sum_{j=0}^{k-1} \, q^j \, ({t_{n-1}}^j g_{n-1}\ldots g_{i+1}g_i){t_n}^{k-j} + \cdots +$ 

\vspace{.1in}
$(q-1)q^{(n-(i+1))k} \, \sum_{j=0}^{k-1} \, q^j \, ({t_i}^jg_i) g_ng_{n-1} \ldots 
g_{i+2}{t_{i+1}}^{k-j} + $

\vspace{.1in}
$q^{(n-i)k}(q-1) \, \sum_{j=0}^{k-1} \, q^j \, ({t_{i-1}}^j)g_ng_{n-1} \ldots
g_{i+1}{t_i}^{k-j}+$

\vspace{.1in}
$q^{(n-i)k}q^k \, g_ng_{n-1} \ldots g_i{t_{i-1}}^k.$

\end{pf}
 
\begin{th}
Every element of ${\cal H}_{n+1}(q,\infty)$ respectively  ${\cal H}_n(q,d)$ is a unique linear
combination of words, each of one of the following types:

\vspace{.1in}
1) $w_{n-1}$

\vspace{.1in}
2) $w_{n-1}g_ng_{n-1} \ldots g_i$

\vspace{.1in}
3) $w_{n-1}g_ng_{n-1} \ldots g_i{t_{i-1}}^k, \  k\in \ZZ$ respectively $k\in \ZZ_d$

\vspace{.1in}
4) $w_{n-1}{t_n}^k, \ k\in \ZZ$ respectively $k\in \ZZ_d$

\vspace{.1in}
\noindent  where $w_{n-1}$ is some word in $\HH$ respectively ${\cal H}_n(q,d)$.  Thus, the above
words furnish an inductive basis for ${\cal H}_{n+1}(q,\infty)$ respectively ${\cal H}_n(q,d)$. 
\end{th} 

\begin{pf}
By Theorem 1 it suffices to show that every element $v\cdot \sigma_n \in \Sigma_2$, where  $v$ is
a product of $t_i$'s and $\sigma_n \in {\cal H}_{n+1}(q)$, can be expressed uniquely in terms of
1), 2), 3) and 4). We prove this by induction on $n$: For $n=0$ there are no $g_i$'s in the word,
so $v\cdot \sigma_0 = t^k\cdot 1$, a word of type 1). Suppose the assertion holds for all basic
words in  $\Sigma_2$ with indices up to $n-1$, and let $w\in \Sigma_2$ such that $w$ contains
elements of index $n$. We examine the different cases:

\vspace{.1in}
\noindent $ \bullet \ \  \ w={t_{i_1}}^{k_1}{t_{i_2}}^{k_2}\ldots
{t_{i_r}}^{k_r}\underline{{t_n}^k \cdot \sigma_{n-1}}, \  
 1\leq i_1<\ldots <i_r< n \mbox{ and }\sigma_{n-1}\in {\cal H}_n(q)$. 

\vspace{.1in}
\noindent Then, by Lemma 1,(v),  
 $w=\underline{{t_{i_1}}^{k_1}\ldots {t_{i_r}}^{k_r}\cdot\sigma_{n-1}}\cdot{t_n}^k =
w_{n-1}{t_n}^k,$ a word of type  

\vspace{.1in}
\noindent  4).

\vspace{.1in}
\noindent $ \bullet \ \ \ w={t_{i_1}}^{k_1}{t_{i_2}}^{k_2}\ldots {t_{i_r}}^{k_r}\cdot \sigma_n,$ \ 
where $i_r<n\mbox{ and }\sigma_n = \sigma_{n-1}\cdot (g_ng_{n-1} \ldots g_i)  \in $

\vspace{.1in}
\noindent $ {\cal H}_{n+1}(q)$. Then $w = \underline{{t_{i_1}}^{k_1}{t_{i_2}}^{k_2}\ldots
{t_{i_r}}^{k_r}\cdot \sigma_{n-1}}\cdot (g_ng_{n-1} \ldots g_i) = w_{n-1}g_ng_{n-1} \ldots g_i,$

\vspace{.1in}  
\noindent  a word of type 2).

\vspace{.1in}
\noindent $\bullet$ \ \ \ Finally, let $w={t_{i_1}}^{k_1}{t_{i_2}}^{k_2}\ldots
{t_{i_r}}^{k_r}{t_n}^k\cdot \sigma_n,$ \ where $ \sigma_n = \sigma_{n-1}\cdot (g_ng_{n-1} \ldots
g_i) \in$ 

\vspace{.1in} 
\noindent $ {\cal H}_{n+1}(q)$.
Then $w={t_{i_1}}^{k_1}\ldots{t_{i_r}}^{k_r}\underline{{t_n}^k\cdot \sigma_{n-1}}\cdot g_ng_{n-1}
\ldots g_i$ $\stackrel{Lemma \, 1,(v)}{=}$  

\vspace{.1in} 
\noindent ${t_{i_1}}^{k_1}\ldots{t_{i_r}}^{k_r}\cdot \sigma_{n-1}\cdot 
\underline{{t_n}^kg_ng_{n-1} \ldots g_i}\stackrel{Lemma \, 2,(ii)}{=} $ 

\vspace{.1in} 
\noindent $w_{n-1}{t_n}^{k-j} + \Sigma w_{n-1}g_ng_{n-1} \ldots g_s{t_{s-1}}^{k-j},$ \ for $j=0,
\ldots, k-1.$

\vspace{.1in} 
\noindent I.e. $w$ is a sum of words of type 4) and type 3). The
uniqueness of these expressions follows from Lemma 1 and Lemma 2. 
\end{pf}

Theorem 2 rephrased weaker says that the elements of the inductive basis contain {\it either}
$g_n$ {\it or} \ ${t_n}^k$ {\it at most once}.
But, as explained in the beginning, our aim is to find an inductive basis for 
${\cal H}_{n+1}$ using the elements 
$t'_i=g_ig_{i-1} \ldots g_1t{g_1}^{-1}\ldots {g_{i-1}}^{-1}{g_i}^{-1},$ 
 as these are the right ones for constructing  Markov traces on $\bigcup_{n=1}^{\infty}{\cal
H}_n.$ We go from the $t_i$'s to the $t'_i$'s via the `intermediate' elements

\[ T_i^k := g_ig_{i-1} \ldots g_1t^k g_1\ldots g_{i-1}g_i, \ k\in \ZZ.\]

\begin{th}
Every element of ${\cal H}_{n+1}(q,\infty)$ respectively ${\cal H}_n(q,d)$ is a unique linear 
combination of words, each of one of the following types:

\vspace{.1in}
1$'$) $w_{n-1}$

\vspace{.1in}
2$'$) $w_{n-1}g_ng_{n-1} \ldots g_i$

\vspace{.1in}
3$'$) $w_{n-1}g_ng_{n-1} \ldots g_iT_{i-1}^k, \  k\in \ZZ$ respectively $k\in \ZZ_d$

\vspace{.1in}
4$'$) $w_{n-1}T_n^k, \  k\in \ZZ$ respectively $k\in \ZZ_d$

\vspace{.1in}
\noindent where $w_{n-1}$ is some word in $\HH$ respectively ${\cal H}_n(q,d)$. 
\end{th}

\begin{pf}
It suffices to show that elements of the inductive basis given in Theorem 2 can be
expressed uniquely as sums of the above words. For this we need the following three lemmas.

\begin{lem} For $k\in \NN$ respectively $k\in \ZZ_{d-1}$ and $\epsilon  \in \{ {\pm 1}\}$ the
following hold  in ${\cal H}_{n+1}(q,\infty)$ respectively ${\cal H}_{n+1}(q,d)$:  

\bigbreak
\noindent ${t_n}^{\epsilon (k+1)} \: = \: \sum_{r_1, \ldots, r_k=0}^{n}
{(q^\epsilon -1)}^{\epsilon_{r_1}+\cdots+\epsilon_{r_k}}  q^{\epsilon (r_1+\cdots+r_k)}\cdot$

\vspace{.1in}
\noindent ${g_n}^\epsilon {g_{n-1}}^\epsilon \ldots {g_1}^\epsilon t^\epsilon ({g_1}^\epsilon
\ldots {g_{n-r_1}}^\epsilon \ldots {g_1}^\epsilon)t^\epsilon \ldots t^\epsilon ({g_1}^\epsilon
\ldots {g_{n-r_k}}^\epsilon \ldots {g_1}^\epsilon)t^\epsilon {g_1}^\epsilon \ldots
{g_{n-1}}^\epsilon {g_n}^\epsilon,$ 

\bigbreak
\noindent where  $\epsilon_{r_i}=1 \mbox{ if } r_i=0, \ldots, n-1,
\, \ \epsilon_n=0 \ \mbox{ and } {g_0}^\epsilon:=1.$

\end{lem}

\begin{pf}
We show the case $\epsilon = +1$ by induction on $k$.  The proof for $\epsilon = -1$ is
completely analogous. For $k=1$ we have: 

\vspace{.1in}
\noindent ${t_n}^2 = g_ng_{n-1} \ldots g_1t \underline{g_1\ldots g_{n-1}g_ng_ng_{n-1} \ldots g_1}  
tg_1\ldots g_{n-1}g_n\stackrel{Lemma \, 1,(iv)}{=}$

\vspace{.1in}
$\sum_{r=0}^{n} {(q-1)}^{\epsilon_r}q^r \, g_ng_{n-1} \ldots g_1t(g_1\ldots g_{n-r} \ldots
g_1)tg_1\ldots g_{n-1} g_n.$ 

\vspace{.1in}
\noindent Assume that the statement holds for any $k\in \NN$. Then for $k+1$ we have:
 
\vspace{.1in}
\noindent ${t_n}^{k+1} = {t_n}^k t_n  \stackrel{by \, induction}{=} \sum_{r_1, \ldots,
r_{k-1}=0}^{n} {(q-1)}^{\epsilon_{r_1}+\cdots+\epsilon_{r_{k-1}}}  q^{r_1+\cdots +r_{k-1}} \, 
g_n\ldots g_1 \cdot$
 
\vspace{.1in}
$t(g_1 \ldots g_{n-r_1}\ldots g_1)t\ldots t(g_1 \ldots g_{n-r_{k-1}}\ldots
g_1)t \underline{g_1\ldots g_n(g_n\ldots g_1} tg_1\ldots g_n)$ 
 
\vspace{.1in}
$ \stackrel{Lemma \, 1,(iv)}{=}  \sum_{r_1, \ldots, r_k=0}^{n}
{(q-1)}^{\epsilon_{r_1}+\cdots+\epsilon_{r_k}}  q^{r_1+\cdots+r_k} \, g_n \ldots g_1\cdot$
 
\vspace{.1in}
$t(g_1 \ldots g_{n-r_1}\ldots g_1)t\ldots t(g_1 \ldots g_{n-r_k}\ldots
g_1)tg_1\ldots g_n.$
\end{pf}

\begin{lem}
For $k\in \NN$  and $\epsilon \in \{ {\pm 1}\}$ the following hold 
 in $\HH$ respectively ${\cal H}_n(q,d)$: 

\vspace{.1in}
$ \begin{array}{rlll}
 (i) & \, t^\epsilon {g_1}^\epsilon t^{\epsilon k}{g_1}^\epsilon &  = & {g_1}^\epsilon
t^{\epsilon k}
 {g_1}^\epsilon t^\epsilon  +  (q^\epsilon -1) t^\epsilon {g_1}^\epsilon t^{\epsilon k}  
+  (1-q^\epsilon) t^{\epsilon k}{g_1}^\epsilon t^\epsilon \ \mbox{ and }  \\ [.1in] 

(ii) & t^{-\epsilon} {g_1}^\epsilon t^{\epsilon k}{g_1}^\epsilon & = & {g_1}^\epsilon 
t^{\epsilon k} {g_1}^\epsilon t^{-\epsilon}  +  (q^\epsilon -1) t^{\epsilon (k-1)}
{g_1}^\epsilon +  (1-q^\epsilon) {g_1}^\epsilon t^{\epsilon (k-1)}.
\end{array} $
\end{lem}

\begin{pf}
We only prove (i) for the case $\epsilon = +1$, by induction on $k$. All other statements are
proved similarly.  For $k=1$ we have   $tg_1tg_1 = g_1tg_1t$. Assume the assertion is correct for
$k$. Then for $k+1$ we have:
 
\vspace{.1in}
\noindent $tg_1t^{k+1}g_1 = tg_1t^kg_1 \underline{g_1^{-1}}tg_1\stackrel{Lemma \, 1,(i)}{=}$
 
\vspace{.1in}
$q^{-1} \, \underline{tg_1t^kg_1} g_1tg_1 + (q^{-1}-1) \, \underline{tg_1t^kg_1}tg_1
\stackrel{induction \, step}{=}$
 
\vspace{.1in}
$q^{-1} \, g_1t^k \underline{g_1tg_1t} g_1 + q^{-1}(q - 1) \, \underline{tg_1t^kg_1} tg_1 + 
q^{-1}(1 - q) \, t^k \underline{g_1tg_1t} g_1 +$
 
\vspace{.1in}
$(q^{-1}-1) \, g_1t^kg_1t^2g_1 +  (q^{-1}-1)(q-1) \, tg_1t^{k+1}g_1 + $
 
\vspace{.1in}
$(q^{-1}-1)(1-q) \, t^kg_1t^2g_1  \stackrel{ rels., \, induction \, step}{=} q^{-1} \,
g_1t^{k+1}g_1t \underline{g_1^2}+ $ 
 
\vspace{.1in}
$ (1- q^{-1}) \, g_1t^kg_1t^2g_1 + (1- q^{-1})(q-1) \, tg_1t^{k+1}g_1 + (1- q^{-1})(1- q) \,
t^kg_1t^2g_1 +$ 

\vspace{.1in}
$(q^{-1}-1) \, t^{k+1}g_1t \underline{g_1^2} + (q^{-1}-1) \, g_1t^kg_1t^2g_1 + 
(q^{-1}-1)(q-1) \, tg_1t^{k+1}g_1 +$

\vspace{.1in}
$(q^{-1}-1)(1-q) \, t^kg_1t^2g_1 \stackrel{Lemma \, 1,(i)}{=}$
 
\vspace{.1in}
$q^{-1}(q-1) \, g_1 \underline{t^{k+1}g_1t g_1} + g_1t^{k+1}g_1t +
(q^{-1}-1)(q-1) \, t^{k+1}g_1t g_1 +$

\vspace{.1in}
$ (q^{-1}-1)q \, t^{k+1}g_1t \stackrel{Lemma \, 1,(v)}{=}$ 
 
\vspace{.1in}
$(1-q^{-1}) \, \underline{g_1^2}t g_1t^{k+1} + g_1t^{k+1}g_1t +
(q^{-1}-1)(q-1) \, t^{k+1}g_1t g_1 + (1-q) \, t^{k+1}g_1t =$ 
 
\vspace{.1in}
$(1-q^{-1})(q-1) \, g_1t g_1t^{k+1} + (1-q^{-1})q \, t g_1t^{k+1} + g_1t^{k+1}g_1t + $
 
\vspace{.1in}
$(q^{-1}-1)(q-1) \, t^{k+1}g_1t g_1 + (1-q) \, t^{k+1}g_1t \stackrel{Lemma \, 1,(v)}{=}$ 
 
\vspace{.1in}
$g_1t^{k+1}g_1t +  (q-1) \, t g_1t^{k+1} +  (1-q) \, t^{k+1}g_1t.$ 
\end{pf}

\begin{lem} [Fundamental Lemma (F.L.)] 
For $i, k\in \NN$  and for $\epsilon \in \{ \pm 1 \}$ the following hold
 in $\HH$ respectively ${\cal H}_n(q,d)$:

\bigbreak
\noindent \ (i) $t^{\epsilon i}{g_1}^\epsilon t^{\epsilon k}{g_1}^\epsilon  =   {g_1}^\epsilon
t^{\epsilon k}{g_1}^\epsilon t^{\epsilon i} + $

\vspace{.1in}
$(q^\epsilon -1) \, [t^\epsilon {g_1}^\epsilon t^{\epsilon
(k+i-1)} + t^{2\epsilon}{g_1}^\epsilon t^{\epsilon (k+i-2)}+ \cdots + t^{\epsilon i}
{g_1}^\epsilon t^{\epsilon k}] +$

\vspace{.1in}
 $(1-q^\epsilon) \, [t^{\epsilon k}{g_1}^\epsilon t^{\epsilon i} +
t^{\epsilon (k+1)}{g_1}^\epsilon t^{\epsilon (i-1)}+ \cdots + t^{\epsilon (k+i-1)}
{g_1}^\epsilon t^\epsilon ] \mbox{ and }$

\bigbreak
\noindent \ (ii) $t^{-\epsilon i}{g_1}^\epsilon t^{\epsilon k}{g_1}^\epsilon  =   {g_1}^\epsilon
t^{\epsilon k}{g_1}^\epsilon t^{-\epsilon i} + $
 
\vspace{.1in}
$(q^\epsilon -1) \, [t^{\epsilon (k-1)} {g_1}^\epsilon t^{-\epsilon (i-1)} +
t^{\epsilon (k-2)}{g_1}^\epsilon t^{-\epsilon (i-2)}+ \cdots + t^{\epsilon (k-i)}
{g_1}^\epsilon ] +$
 
\vspace{.1in}
 $(1-q^\epsilon) \, [t^{-\epsilon (i-1)}{g_1}^\epsilon t^{\epsilon (k-1)} +
t^{-\epsilon (i-2)}{g_1}^\epsilon t^{\epsilon (k-2)}+ \cdots + {g_1}^\epsilon t^{\epsilon (k-i)}
].$

\end{lem}

\begin{pf}
We prove (i) for the case $\epsilon = +1$, by induction on $i$.  The proof for $\epsilon = -1$ is
completely analogous.  For $i=1$ the assertion
is true by Lemma 4,(i). Assume it holds for $i$. Then for $i+1$ we have: 
 
\vspace{.1in}
\noindent $t^{i+1}g_1t^kg_1 = t\underline{t^ig_1t^kg_1} \stackrel{induction \, step}{=} 
\underline{tg_1t^kg_1}t^i + $
 
\vspace{.1in}
$(q-1) \, [t^2g_1t^{k+i-1} + t^3g_1t^{k+i-2}+ \cdots + t^{i+1}g_1t^k] + $

\vspace{.1in}
$(1-q) \, [t^{k+1}g_1t^i + t^{k+2}g_1t^{i-1}+ \cdots + t^{k+i}g_1t] \stackrel{Lemma \, 4,(i)}{=}$

\vspace{.1in}
$g_1t^kg_1t^{i+1} + (q-1) \, tg_1t^{k+i} + (1-q) \, t^kg_1t^{i+1} +$

\vspace{.1in}
$(q-1) \, [t^2g_1t^{k+i-1} + t^3g_1t^{k+i-2}+ \cdots + t^{i+1}g_1t^k] + $

\vspace{.1in}
$(1-q) \, [t^{k+1}g_1t^i + t^{k+2}g_1t^{i-1}+ \cdots + t^{k+i}g_1t].$

\end{pf}

We go back now to the proof of Theorem 3. By Lemma 3, a typical summand of
${t_n}^{\epsilon(k+1)} \in {\cal H}_{n+1}$ is: 

\bigbreak
\noindent ${g_n}^\epsilon\ldots {g_1}^\epsilon t^{\epsilon\lambda_1}
({g_1}^\epsilon \ldots {g_{n-l_1}}^\epsilon\ldots {g_1}^\epsilon)t^{\epsilon\lambda_2}\ldots
t^{\epsilon\lambda_N} ({g_1}^\epsilon \ldots {g_{n-l_N}}^\epsilon\ldots
{g_1}^\epsilon)t^{\epsilon\lambda_{N+1}} {g_1}^\epsilon\ldots {g_n}^\epsilon,$

\bigbreak
\noindent $ \mbox{ where } \lambda_1,\lambda_2, \ldots, \lambda_{N+1} \in
\NN$ such that $ \lambda_1+\cdots+ \lambda_{N+1} = k+1 \mbox{ and } l_i<n \mbox{ for } i=1,
\ldots, N$ (since the cases $l_i=n$ are incorporated in $t^{\epsilon\lambda_i}$). \ In order to
prove the theorem we want to show that such a word can be expressed in terms of words of the form
1$'$), 2$'$), 3$'$) and 4$'$). This is a very  slow process as we shall readily see. In order to
obtain an inductive argument on the number $N+1$ of the intermediate powers of $t$, we  show
first the following, seemingly more general result, where an unsymmetric expression appears
also in the word. It is  'seemingly more general' because this unsymmetry of the word appears
anyhow in a later stage of the calculations. 

\begin{propn}
Let  $ k\in \NN$ respectively $ k\in \ZZ_{d-1}$, $\epsilon \in \{ \pm 1 \}, \ l, m, l_2, \ldots,
l_N\leq n$ and let $\lambda_1, \lambda_2, \ldots, \lambda_{N+1} \in \NN$  such that 
$\lambda_1+\cdots+ \lambda_{N+1} = k+1 $.  Then it holds in ${\cal H}_{n+1}(q,\infty)$ 
respectively ${\cal H}_{n+1}(q,d)$ that words of the form:  
 
\vspace{.1in}
\noindent $w_{n-1} {g_n}^\epsilon\ldots {g_1}^\epsilon t^{\epsilon\lambda_1}
({g_1}^\epsilon \ldots {g_l}^\epsilon)
({g_1}^\epsilon \ldots {g_m}^\epsilon\ldots {g_1}^\epsilon)t^{\epsilon\lambda_2} ({g_1}^\epsilon
\ldots {g_{l_2}}^\epsilon\ldots {g_1}^\epsilon)t^{\epsilon\lambda_3}\ldots$

\vspace{.1in}
\noindent $t^{\epsilon\lambda_{N+1}}  {g_1}^\epsilon\ldots {g_n}^\epsilon$
 
\vspace{.1in}
\noindent where {\it only} between the first two powers of \ $t$ appears the unsymmetric
expression  $({g_1}^\epsilon \ldots {g_l}^\epsilon) ({g_1}^\epsilon \ldots {g_m}^\epsilon\ldots
{g_1}^\epsilon),$ can be expressed as sums of words of the form 1$'$), 2$'$), 3$'$) and 4$'$).
Note that if $l=0$ we obtain the generic summand of ${t_n}^{\epsilon(k+1)}$.  
\end{propn}

\begin{pf}
We prove the statement for $\epsilon=+1$ by induction on the number $N+1$ of intermediate
powers of $t$. The proof for  $\epsilon=-1$ is completely analogous.   For $N=0$ we have   
$w_{n-1} g_n\ldots g_1 t^{\lambda_1}g_1\ldots g_n,$ where $\lambda_1=k+1$ i.e.
$w_{n-1}T_n^{k+1}.$  Suppose the assertion holds for $N$. Then for $N+1$ we have:    
 
\vspace{.1in}
\noindent $A=w_{n-1} g_n\ldots g_1 t^\lambda (g_1\ldots g_l)(g_1\ldots g_m\ldots g_1)t^\mu
(g_1\ldots g_{l_2} \ldots g_1)t^{\lambda_3}\ldots t^{\lambda_{N+1}} g_1 \ldots g_n$

\vspace{.1in}
\noindent $=w_{n-1} g_n\ldots g_1 t^\lambda \underline{(g_1\ldots g_l)}(g_m\ldots g_1\ldots
g_m)t^\mu (g_{l_2}\ldots g_1 \ldots g_{l_2})t^{\lambda_3}\ldots t^{\lambda_{N+1}} g_1 \ldots g_n.$

\vspace{.1in}
\noindent Here we also use the symbol `$\sum$' to mean `linear combination of words of the
type',  the symbol `$w_{n-1}$' for {\it not } always the same word in ${\cal H}_n$, and, in order
to shorten the words, we substitute the expression   
$ g_{l_2}\ldots g_1 \ldots g_{l_2}t^{\lambda_3}\ldots t^{\lambda_{N+1}} g_1 \ldots g_n$ by
$S$.   

\vspace{.1in}
We proceed by examining the cases $l<m, \ l>m$ and $ l=m.$

\vspace{.1in}
\noindent $\bullet$ \ For $l<m$ we have: 
 
\vspace{.1in}
\noindent $A= w_{n-1} g_n\ldots g_1 t^\lambda \underline{g_m\ldots g_2}g_1\ldots g_m(g_1\ldots
g_l)t^\mu\cdot S=$ 
 
\vspace{.1in}
$w_{n-1}(g_{m-1}\ldots g_1) g_n\ldots g_1 t^\lambda g_1\ldots g_m \underline{g_1}g_2\ldots 
g_lt^\mu \cdot S \stackrel{m>1}{=}$

\vspace{.1in}
$w_{n-1}(g_{m-1}\ldots g_1g_1) g_n\ldots g_1 t^\lambda g_1\ldots g_m  g_2\ldots g_l
\underline{t^\mu}\cdot S =$
 
\vspace{.1in}
$\underline{w_{n-1}(g_{m-1}\ldots g_1^2)} g_n\ldots \underline{g_1 t^\lambda g_1t^\mu}   
g_2\ldots g_m  g_2\ldots g_l\cdot S \stackrel{F.L.}{=} $
 
\vspace{.1in}
$w_{n-1}g_n\ldots g_2 \underline{t^\mu}g_1 t^\lambda g_1g_2\ldots g_mg_2\ldots g_l\cdot S +$
 
\vspace{.1in}
$\sum_{i+j=\lambda+\mu}^{} w_{n-1}g_n\ldots g_2 \underline{t^i}g_1 \underline{t^j}
g_2\ldots g_mg_2\ldots g_l\cdot S =$
 
\vspace{.1in}
$w_{n-1}t^\mu g_n\ldots g_1 t^\lambda g_1\ldots g_m \underline{g_2\ldots g_l}\cdot S +$
 
\vspace{.1in}
$\sum_{i+j=\lambda+\mu}^{} w_{n-1}t^i g_n\ldots g_1 \underline{g_2\ldots g_mg_2\ldots g_l}t^j 
\cdot S\stackrel{Lemma \, 1,(ii), \, l<m }{=}$

\vspace{.1in}
$\underline{(w_{n-1}t^\mu g_2\ldots g_l)}g_n\ldots g_1 t^\lambda g_1\ldots g_m \cdot S+$
 
\vspace{.1in}
$\sum_{i+j=\lambda+\mu}^{} \underline{(w_{n-1}t^i g_1\ldots g_{m-1}g_1\ldots g_{l-1})} g_n\ldots
g_1 t^j \cdot S =$
 
\vspace{.1in}
$(w_{n-1}g_n\ldots g_1 t^\lambda g_1\ldots g_m\cdot S + \sum_{i+j=\lambda+\mu}^{}
w_{n-1}g_n\ldots g_1 t^j \cdot S$    
 
\vspace{.1in}
\noindent and the number of intermediate powers of $t$ has reduced to $N$ in all summands  of
${t_n}^{k+1}$. 

\vspace{.1in}
\noindent $\bullet$ \ For $l>m$ we have: 
 
\vspace{.1in}
\noindent $A= w_{n-1} g_n\ldots g_1 t^\lambda(g_1\ldots g_l) \underline{g_m\ldots g_1\ldots
g_m}t^\mu \cdot S  \stackrel{m<l, \, Lemma \, 1,(ii)}{=}$
 
\vspace{.1in}
$\underline{(w_{n-1}g_m\ldots g_1\ldots g_m)}g_n\ldots g_1 t^\lambda g_1\ldots g_l
\underline{t^\mu}\cdot S =$
 
\vspace{.1in}
$w_{n-1}g_n\ldots \underline{g_1 t^\lambda g_1t^\mu}g_2\ldots g_l \cdot S 
\stackrel{F.L.}{=} w_{n-1}g_n\ldots g_2\underline{t^\mu} g_1 t^\lambda g_1\ldots g_l\cdot S +$

\vspace{.1in}
$\sum_{i+j=\lambda+\mu}^{} w_{n-1} g_n\ldots g_2 \underline{t^i}g_1 \underline{t^j}
g_2\ldots g_l\cdot S =$
 
\vspace{.1in}
$\underline{w_{n-1}t^\mu} g_n\ldots g_1 t^\lambda g_1\ldots g_l\cdot S
+ \sum_{i+j=\lambda+\mu}^{} w_{n-1}t^i g_n\ldots g_1 \underline{(g_2\ldots g_l)}t^j\cdot S 
\stackrel{Lemma \, 1,(ii)}{=}$
 
\vspace{.1in}
$w_{n-1} g_n\ldots g_1 t^\lambda g_1\ldots g_l\cdot S + \sum_{i+j=\lambda+\mu}^{}
\underline{(w_{n-1}t^i g_1\ldots g_{l-1})}g_n\ldots g_1 t^j\cdot S =$

\vspace{.1in}
$w_{n-1} g_n\ldots g_1 t^\lambda g_1\ldots g_l + \sum_{i+j=\lambda+\mu}^{} w_{n-1} g_n\ldots g_1
t^j$
 
\vspace{.1in}
\noindent and the number of intermediate powers of $t$ has reduced to $N$ in all summands  of
${t_n}^{k+1}$.

\vspace{.1in}
\noindent $\bullet \ \mbox{ Finally if } l=m \mbox{ we have: }$ 
 
\vspace{.1in}
\noindent $A= w_{n-1} g_n\ldots g_1 t^\lambda \underline{(g_1\ldots g_m) g_m\ldots g_1}\ldots
g_mt^\mu \cdot S  \stackrel{Lemma 1, \, (iv)}{=}$ 
 
\vspace{.1in}
$ w_{n-1} g_n\ldots g_1 t^\lambda \underline{g_2\ldots g_m}t^\mu \cdot S  +$

\vspace{.1in}
$\sum_{r=0}^{m-1} w_{n-1} g_n\ldots g_1 t^\lambda (\underline{g_{m-r}\ldots g_2}
g_1\ldots g_{m-r})  g_2\ldots g_m\underline{t^\mu} \cdot S  \stackrel{Lemma 1, \, (ii)}{=} $
 
\vspace{.1in}
$\underline{(w_{n-1}g_1\ldots g_{m-1})} g_n\ldots g_1 t^{\lambda +\mu}\cdot S +$

\vspace{.1in}
$\sum_{r=0}^{m-1} \underline{(w_{n-1}g_{m-r-1}\ldots g_1) } g_n\ldots 
\underline{g_1 t^\lambda g_1t^\mu}(g_2\ldots g_{m-r}) g_2\ldots g_m \cdot S \stackrel{F.L.}{=} $
 
\vspace{.1in}
$ w_{n-1} g_n\ldots g_1 t^{\lambda +\mu}\cdot S  + \sum_{r=0}^{m-1}w_{n-1}g_n\ldots
g_2\underline{t^\mu} g_1 t^\lambda (g_1\ldots g_{m-r}) g_2\ldots g_m \cdot S + $
 
\vspace{.1in}
$\sum_{i+j=\lambda+\mu}^{}\sum_{r=0}^{m-1} w_{n-1}g_n\ldots
g_2\underline{t^i} g_1 \underline{t^j} (g_2\ldots g_{m-r}) g_2\ldots g_m \cdot S = $
  
\vspace{.1in}
$ w_{n-1} g_n\ldots g_1 t^{\lambda +\mu} \cdot S + \sum_{r=0}^{m-1}w_{n-1}t^\mu g_n\ldots
g_1t^\lambda (g_1\ldots g_{m-r})  \underline{g_2\ldots g_{m-r-1}}\ldots g_m \cdot  $
 
\vspace{.1in}
$S + \sum_{i+j=\lambda+\mu}^{}\sum_{r=0}^{m-1} w_{n-1}t^i g_n\ldots g_1            
\underline{(g_2\ldots g_{m-r})g_2\ldots g_m}  t^j \cdot S \stackrel{Lemma 1, \, (ii)}{=} $
 
\vspace{.1in}
$ w_{n-1} g_n\ldots g_1 t^{\lambda +\mu} \cdot S +$

\vspace{.1in}
$\sum_{r=0}^{m-1}\underline{(w_{n-1}t^\mu g_2\ldots g_{m-r-1})} g_n\ldots g_1t^\lambda  
(g_1\ldots g_{m-r-1} \underline{g_{m-r}^2}g_{m-r+1} \ldots g_m) \cdot S + $
 
\vspace{.1in}
$\sum_{i+j=\lambda+\mu}^{}\sum_{r=0}^{m-1} \underline{(w_{n-1}t^ig_1\ldots g_{m-r-1}g_1\ldots
g_{m-1})} g_n\ldots g_1 t^j \cdot S = $

\vspace{.1in}
$ w_{n-1} g_n\ldots g_1 t^{\lambda +\mu} \cdot S +  \sum_{r=0}^{m-1}w_{n-1}g_n\ldots
g_1t^\lambda      (g_1\ldots g_{m-r-1}\underline{g_{m-r+1}\ldots g_m}) \cdot S + $    

\vspace{.1in}
$ \sum_{r=0}^{m-1}w_{n-1}g_n\ldots g_1t^\lambda (g_1\ldots g_m) \cdot S + 
\sum_{i+j=\lambda+\mu}^{}\sum_{r=0}^{m-1} w_{n-1} g_n\ldots g_1 t^j \cdot S = $
 
\vspace{.1in}
$ w_{n-1} g_n\ldots g_1 t^{\lambda +\mu} \cdot S + \sum_{r=0}^{m-1}
\underline{w_{n-1}g_{m-r+1}\ldots g_m}g_n\ldots g_1t^\lambda (g_1\ldots g_{m-r-1}) \cdot
S+ $        

\vspace{.1in}
$\sum_{r=0}^{m-1}w_{n-1}g_n\ldots g_1t^\lambda (g_1\ldots g_m) \cdot S +
\sum_{i+j=\lambda+\mu}^{}\sum_{r=0}^{m-1} w_{n-1} g_n\ldots g_1 t^j \cdot S = $

\vspace{.1in}
$ w_{n-1} g_n\ldots g_1 t^{\lambda +\mu} \cdot S + \sum_{r=0}^{m-1}
w_{n-1}g_n\ldots g_1t^\lambda (g_1\ldots g_{m-r-1}) \cdot
S+ $        

\vspace{.1in}
$\sum_{r=0}^{m-1}w_{n-1}g_n\ldots g_1t^\lambda (g_1\ldots g_m) \cdot S +
\sum_{i+j=\lambda+\mu}^{}\sum_{r=0}^{m-1} w_{n-1} g_n\ldots g_1 t^j \cdot S $
 
\vspace{.1in}
\noindent and the number of the intermediate powers of $t$ has reduced to $N$ in all summands of
${t_n}^{k+1}$.

\end{pf}

\noindent We can now conclude the proof of Theorem 3, since for the different possibilities of a
word $w \in {\cal H}_{n+1}$ we have: 
 
\begin{itemize}
\item[Case 1.] $ \ \mbox{ If } w= w_{n-1} \mbox{ or } w= w_{n-1} g_n\ldots g_i \mbox{ for } i=0,
\ldots, n$  there is nothing to show.  
\item[Case 2.] $ \ \mbox{ If } w= w_{n-1}t_n^k, \ k\in\ZZ$ respectively  $\ZZ_d$, then by
Proposition 3, $w$ is a unique linear combination of words of type 1$'$), 2$'$), 3$'$) and 4$'$). 
\item[Case 3.] $ \ \mbox{ Finally, if } w= w_{n-1} g_n\ldots g_{i+1}t_i^k, \ 
k\in\ZZ$ respectively $\ZZ_d$, by Proposition 3,   $t_i^k$ is written in terms of words \ 
 $w_{i-1}, \, w_{i-1}g_i\ldots g_r$ \ for  $r\leq i, \ w_{i-1}g_i\ldots g_{r+1}T_r^k$  
 \  and  \ $w_{i-1}T_i^k.$ \ Therefore $w$ can be written uniquely in terms of the words 
 
\vspace{.1in}
 $w_{n-1}g_n\ldots g_{i+1}w_{i-1}g_i\ldots g_r\mbox{ for } r=0, \ldots, i$,
 
\vspace{.1in}
 $w_{n-1}g_n\ldots g_{i+1}w_{i-1}g_i\ldots g_{r+1}T_r^k$ \ and 
 
\vspace{.1in}
 $w_{n-1}g_n\ldots g_{i+1}w_{i-1}T_i^k.$ 
 
\vspace{.1in}
\noindent $w_{i-1}$ commutes with $g_n\ldots g_{i+1}$, unless $i=0$, where the word is already
arranged in a trivial manner. So the above words reduce to the types 
$w_{n-1}g_n\ldots g_r$ or $w_{n-1}g_n\ldots g_{j+1}T_j^k.$
\end{itemize} 
\end{pf}

Theorem 3 rephrased weaker says that the elements of the inductive basis contain {\it either}
$g_n$ {\it or} \ $T_n^k$ {\it at most once}.
We can now pass easily to the inductive basis that we need for constructing Markov
traces on  $\bigcup_{n=1}^{\infty}{\cal H}_n$. Indeed we have the following: 

\begin{th}
Every element of ${\cal H}_{n+1}(q,\infty)$ respectively  ${\cal H}_{n+1}(q,d)$  can be written
uniquely as a linear combination of words, each of one of the following types:

\vspace{.1in}
1$''$) $w_{n-1}$

\vspace{.1in}
2$''$) $w_{n-1}g_ng_{n-1} \ldots g_i$

\vspace{.1in}
3$''$) $w_{n-1}g_ng_{n-1} \ldots g_{i+1}{t'_i}^k, k\in \ZZ$ respectively $\ZZ_d$

\vspace{.1in}
4$''$) $w_{n-1}{t'_n}^k, k\in \ZZ$ respectively $\ZZ_d$

\vspace{.1in}
where $w_{n-1}$ is some word in $\HH$ respectively ${\cal H}_n(q,d)$. 

\end{th}  

\begin{pf} 
By Theorem 3 it suffices to show that expressions of the forms 3$'$) and 4$'$) can be written
(uniquely) in terms of 1$''$), 2$''$), 3$''$) and 4$''$). Indeed, for $k\in \ZZ$,
let 

\vspace{.1in}
\noindent $w = w_{n-1}g_ng_{n-1} \ldots g_{i+1}T_i^k  = w_{n-1}g_ng_{n-1} \ldots g_{i+1}g_i\ldots
g_1 t^k \underline{g_1\ldots g_i}.$ 

\vspace{.1in}
\noindent We apply the relation $g_r = q\cdot g_r^{-1} + (q-1)\cdot 1$ to all letters of the word 
$g_1\ldots g_i$ to get:  

\vspace{.1in}
\noindent $w = w_{n-1}g_n \ldots g_{i+1}g_i\ldots g_1 t^k g_1^{-1}\ldots g_i^{-1} +  \sum_{}^{} 
w_{n-1}g_n \ldots g_1 t^k g_{j_1}^{-1}\ldots g_{j_k}^{-1},$ 

\vspace{.1in}
\noindent where in the words 
$g_{j_1}^{-1}\ldots g_{j_k}^{-1}$ there are possible gaps of indices. Let the closest to
$t^k$ gap occur at the index $\rho$; then   

\vspace{.1in}
\noindent $w = w_{n-1}g_n \ldots g_{i+1} {t'_i}^k +  \sum_{}^{}  w_{n-1}g_n \ldots g_1 t^k
g_1^{-1}\ldots g_{\rho-1}^{-1}   \underline{g_{\rho+1}^{-1}\ldots g_{j_k}^{-1}}=$ 

\vspace{.1in}
$ w_{n-1}g_n \ldots g_{i+1} {t'_i}^k + \sum_{}^{} \underline{(w_{n-1}g_\rho^{-1}\ldots
g_{j_k-1}^{-1})} g_n \ldots g_1 t^k g_1^{-1}\ldots g_{\rho-1}^{-1} =$ 

\vspace{.1in}
$ w_{n-1}g_n \ldots g_{i+1} {t'_i}^k + \sum_{}^{} w_{n-1} g_n \ldots g_\rho {t'_{\rho-1}}^k.$

\vspace{.1in}
\noindent Hence $w$ is a sum of words of type 3$''$. In the case where  $w = w_{n-1}T_n^k, \ k\in
\ZZ$, we apply the same reasoning as above. 

\end{pf}

 Theorem 4 rephrased weaker says that the elements of the inductive basis
contain {\it either} $g_n$ {\it or} \ ${t'_n}^k$ {\it at most once}. Notice also that if we were
working on the level of the Iwahori-Hecke algebra ${\cal H}_n(q,Q)$, we would omit Theorem 3. 

\begin{rem} \rm All three inductive bases of ${\cal H}_{n+1}(q,\infty)$ respectively ${\cal
H}_{n+1}(q,d)$ given in Theorems 2, 3 and 4 induce the same complete set of right coset
representatives, $S_{n+1}$, of $W_{n,\infty}$ respectively $W_{n,d}$ in $W_{n+1,\infty}$ 
respectively $W_{n+1,d}$, namely:

\[ \begin{array}{c} 
 S_{n+1} := \{s_ns_{n-1} \ldots s_i \, | \, i=1, \ldots, n\}\bigcup \\ [.1in] 

 \{ s_ns_{n-1} \ldots s_1 t^k s_1\ldots s_i \, | \, i=1, \ldots, n-1,  \, k\in \ZZ \mbox{
respectively } k\in  \ZZ_d, k \neq 0\} \bigcup \\ [.1in] 

 \{ {t_n}^k \, | \, k\in \ZZ \mbox{ respectively } k\in  \ZZ_d \}. 
 \end{array} \]
\end{rem}

We now give the final result of this section, namely, a basic set of ${\cal H}_{n+1}$
which is a proper subset of $\Sigma_1$. 

\begin{th}

The set
 \[\Sigma = \{ {t'_{i_1}}^{k_1}{t'_{i_2}}^{k_2}\ldots {t'_{i_r}}^{k_r}\cdot \sigma \}\]
for   $1\leq i_1<i_2<\ldots <i_r\leq n, \ k_1, \ldots, k_r \in \ZZ$  respectively $\ZZ_d$ and 
$\sigma \in {\cal H}_{n+1}(q)$ forms a basis in ${\cal H}_{n+1}(q,\infty)$ respectively ${\cal
H}_{n+1}(q,d)$.  
\end{th}

\begin{pf} 
By Theorem 4 it suffices to show that words in the inductive basis 1$''$), 2$''$), 3$''$) and
4$''$)
 can be written in terms of elements of $\Sigma$. Indeed, by induction on $n$ we have: if $n=0$
the only non-empty words are powers of $t$, which are of type 4$''$) and which are  elements of
$\Sigma$ trivially. Assume the result holds for $n-1$. Then for $n$ we have: 

\begin{itemize}
\item[Case 1.]  \ If $w = w_{n-1}$ there is nothing to show (by induction).

\item[Case 2.]  \  If $w = w_{n-1}g_n\ldots g_i$, then, by induction $ w_{n-1} =
{t'_{i_1}}^{k_1}\ldots {t'_{i_r}}^{k_r}\cdot \sigma$, a word of $\Sigma$ restricted on 
${\cal H}_n.$  Thus $ w = {t'_{i_1}}^{k_1}\ldots {t'_{i_r}}^{k_r}\cdot \sigma \cdot g_n\ldots g_i
\in \Sigma ,$ since $\sigma \cdot g_n\ldots g_i$ is an element of the canonical basis of 
${\cal H}_{n+1}(q)$. 

\item[Case 3.]  \  If $w = w_{n-1}g_n\ldots g_{i+1}{t'_{i}}^k$, then, by induction step 
$w_{n-1} = {t'_{i_1}}^{k_1}\ldots {t'_{i_r}}^{k_r}\cdot \sigma$, a word of $\Sigma$ restricted on 
${\cal H}_n,$ so 

\vspace{.1in}
$ w = {t'_{i_1}}^{k_1}\ldots {t'_{i_r}}^{k_r}\cdot \sigma\cdot \underline{g_n\ldots
g_{i+1}{t'_{i}}^k} \stackrel{Lemma 1, \, (vi)}{=}  $

\vspace{.1in}
${t'_{i_1}}^{k_1}\ldots {t'_{i_r}}^{k_r}\cdot
\underline{\sigma\cdot {t'_{n}}^k} g_n\ldots g_{i+1} \stackrel{Lemma 1, \, (vi)}{=}$

\vspace{.1in}
${t'_{i_1}}^{k_1}\ldots{t'_{i_r}}^{k_r}{t'_{n}}^k \cdot \sigma\cdot g_n\ldots g_i$. 

\vspace{.1in}
Now $\sigma\cdot g_n\ldots g_i$ is a basic element of ${\cal H}_{n+1}(q)$, thus $w \in \Sigma$.

\item[Case 4.] \ Finally, if $w = w_{n-1}{t'_{n}}^k,$ by induction step we have 
$ w_{n-1} = {t'_{i_1}}^{k_1}\ldots {t'_{i_r}}^{k_r}\cdot \sigma$, a word of $\Sigma$ restricted
on  ${\cal H}_n$. Then 

\vspace{.1in}
$w = {t'_{i_1}}^{k_1}\ldots {t'_{i_r}}^{k_r}\cdot \underline{\sigma\cdot
{t'_{n}}^k}\stackrel{Lemma 1, \, (vi)}{=}  {t'_{i_1}}^{k_1}\ldots {t'_{i_r}}^{k_r}{t'_{n}}^k\cdot
\sigma \in \Sigma.$ 
\end{itemize}
\end{pf}

\section{Construction of Markov traces}

The aim of this section is to construct Markov linear traces on the generalized and on each
level of the cyclotomic Iwahori-Hecke algebras of ${\cal B}$-type.  
 As these algebras are quotients of the braid groups, the constructed traces will actually attach
to each braid a Laurent polynomial. The traces as well as the strategy of their construction are
based on and include as special case the one constructed on the classical ${\cal B}$-type Hecke
algebras  in \cite{L1}, \cite{L2} (Theorem 5), which in turn was based on Ocneanu's trace on Hecke
algebras of  ${\cal A}$-type, cf. \cite{J} (Theorem 5.1). In the next section we combine these
results with the Markov braid equivalence for knots in a solid torus, so as to obtain analogues
of the homfly-pt polynomial for the solid torus. 
\bigbreak

Let  ${\cal R} = \ZZ \, [q^{\pm 1}, u_1^{\pm 1}, \ldots, u_d^{\pm 1}, \ldots]$ and let ${\cal
H}_n$ denote either  ${\cal H}_n(q,\infty)$ or  ${\cal H}_n(q,d)$. Note  that the natural
inclusion of the group $B_{1,n}$ into $B_{1,n+1}$ (geometrically, by adding one more strand at
the end of the braid) induces  a natural inclusion of  ${\cal H}_n$ into ${\cal H}_{n+1}$. 
Therefore it makes sense to consider ${\cal B} := \bigcup_{n=1}^{\infty}B_{1,n}$ and ${\cal H}
:= \bigcup_{n=1}^{\infty}{\cal H}_n$.  Then we have the following result:

\begin{th} Given $z, s_k$, specified elements in ${\cal R}$ with  $k\in \ZZ$ respectively 
$\ZZ_d$ and  $k\neq 0$, there exists a unique linear trace function
 
\[ tr: \ {\cal H} := \bigcup_{n=1}^{\infty}{\cal H}_n  \longrightarrow  
{\cal R}(z, s_k), \ k\in \ZZ \mbox{ respectively } \ZZ_d \]

determined by the rules:

\[ \begin{array}{lll} 
1) & tr(ab)=tr(ba)              & a,b \in {\cal H}_n \\  
2) & tr(1)=1                    & \mbox{for all }  {\cal H}_n \\ 
3) & tr(ag_n)=z\, tr(a)         & a \in {\cal H}_n \\ 
4) & tr(a{t_n'}^k)=s_k \, tr(a) & a \in {\cal H}_n, \ k\in \ZZ \mbox{ respectively }\ZZ_d  
\end{array} \]

\end{th}

\begin{pf} 
The idea of the proof of Theorem 6 is to construct $tr$ on $\bigcup_{n=1}^{\infty}{\cal H}_n$
 inductively using Theorem 4 and the two last rules of the statement above. For this we need the
following lemma. In order to avoid confusion with the indices we introduce here the symbol
`$Z$' to mean `$\ZZ$' or `$\ZZ_d$' respectively.  

\begin{lem}
The map 

\vspace{.1in}
$ \begin{array}{lcl}
c_n : \ ({\cal H}_n \bigotimes_{{\cal H}_{n-1}} {\cal H}_n) \ \bigoplus_{ k\in Z} {\cal H}_n  &
\longrightarrow  & {\cal H}_{n+1}  \\ [.1in] 
\mbox{given by } \ \ \ \ \ \ \ \ \ \ \ \ \ \ \ \ \ c_n(a \otimes b \oplus_k e_k) & := & ag_nb +
\sum_{k\in Z} e_k{t'_n}^k \end{array} $

\vspace{.1in}
\noindent is an isomorphism of \ $({\cal H}_n,{\cal H}_n)$-bimodules.

\end{lem}

\begin{pf} It follows from Theorem 4 that the set ${\cal L}_n$ below provides  a basis of ${\cal
H}_n$ as a free ${{\cal H}_{n-1}}$-module (compare  with Remark 2 for $W_{n+1}$):

\vspace{.1in}
\noindent ${\cal L}_n := \{g_{n-1}g_{n-2} \ldots g_i \, | \, i=1, \ldots, n-1\} \bigcup \{
{t'_{n-1}}^k \, | \, k\in Z \}\bigcup $

\vspace{.1in}
$ \{ g_{n-1}g_{n-2} \ldots g_1 t^k {g_1}^{_1}\ldots {g_i}^{_1} \, | \,  i=1, \ldots, n-2, \ k\in
Z, \ k \neq 0\}.$

\vspace{.1in}
\noindent Then we have: \ \ ${\cal H}_n  \ = \ \bigoplus_{b\in {\cal L}_n}{\cal H}_{n-1} \cdot
b,$

\vspace{.1in}
\noindent and using the universal property of tensor product we obtain:

\vspace{.1in}
$ \begin{array}{lcl}
{\cal H}_n \bigotimes_{{\cal H}_{n-1}} {\cal H}_n & = & {\cal H}_n \bigotimes_{{\cal H}_{n-1}} 
(\bigoplus_{b\in {\cal L}_n} {\cal H}_{n-1} \cdot b) \\ [.1in] 
 & = & \bigoplus_{b\in {\cal L}_n}({\cal H}_n  \bigotimes_{{\cal H}_{n-1}} {\cal H}_{n-1} \cdot b)
\\ [.1in] 
& = &  \bigoplus_{b\in {\cal L}_n}{\cal H}_n \cdot b.
\end{array} $

\vspace{.1in}
\noindent Therefore:

\vspace{.1in}
\noindent ${\cal H}_n \bigotimes_{{\cal H}_{n-1}} {\cal H}_n \bigoplus_{ k\in Z} {\cal H}_n  \ =
\ \bigoplus_{b\in {\cal L}_n}{\cal H}_n \cdot b \bigoplus_{ k\in Z} {\cal
H}_n.$ 

\vspace{.1in}
\noindent Applying now the same reasoning as above, the set ${\cal L}_{n+1}$ below provides a
basis of ${\cal H}_{n+1}$ as a free ${{\cal H}_n}$-module:

\vspace{.1in}
\noindent ${\cal L}_{n+1} := \{g_ng_{n-1} \ldots g_i \, | \, i=1, \ldots, n\}\bigcup \{ {t'_n}^k
\, | \, k\in Z \}\bigcup $

\vspace{.1in}
$ \{ g_ng_{n-1} \ldots g_1 t^k {g_1}^{-1}\ldots {g_i}^{-1}  \, | \, i=1, \ldots, n-1, \ k\in Z, \
k \neq 0\}.$

\vspace{.1in}
\noindent The latter isomorphism then proves that $c_n$ is indeed an isomorphism of  
 \ $({\cal H}_n,{\cal H}_n)$-bimodules, since it corresponds bijectively basic elements to
 elements of the set ${\cal L}_{n+1}$.  

\end{pf}

We can now define inductively a trace, $tr$, on ${\cal H} = \bigcup_{n=1}^{\infty}{\cal H}_n$ as
follows: assume $tr$ is defined on ${\cal H}_n$ and let $x\in {\cal H}_{n+1}$ be an arbitrary
element. By Lemma 6 there exist $a, b, e_k\in {\cal H}_n, \ k\in Z$, such that  

\[x := c_n( a \otimes b \oplus_k e_k).\] 

\noindent Define now:

\[ tr(x) \ := \ z \cdot tr(ab) + tr(e_0) + \sum_{k\in Z} s_k
\cdot tr(e_k).  \]

\noindent Then $tr$ is well-defined. Furthermore, it satisfies the rules 2), 3) and 4)
of the statement of Theorem 6. Rule 3) reflects the Markov property  (recall the discussion in 
Introduction), and therefore, if the function $tr$ is a trace then it is in particular a Markov
trace. In fact one can check easily using induction and linearity, that $tr$ satisfies
the following seemingly stronger condition:  

\[ (3') \ \ tr(ag_nb)=z\, tr(ab), \ \mbox{ for any } a, b \in {\cal H}_n. \]

In order to prove the existence of $tr$, it remain to prove the conjugation property, i.e.
that  $tr$ is indeed a trace. We show this by examining case by case the different
possibilities. 
Before continuing with the proof, we note that having proved the existence, the uniqueness of $tr$
follows immediately, since for any $x \in {\cal H}_{n+1}, \ tr(x)$ can be clearly computed
inductively using rules 1),  2), 3), 4) and linearity.
 
\bigbreak
\noindent We now proceed with checking that $tr(ax)=tr(xa)$ for all $a, x \in {\cal H}$. Since
$tr$ is defined inductively the assumption holds for all $a, x\in {\cal H}_n$, and we shall
show that $tr(ax)=tr(xa)$ for  $a, x\in {\cal H}_{n+1}$. For this it suffices to consider 
$a \in {\cal H}_{n+1}$ arbitrary and $x$ one of the generators of ${\cal H}_{n+1}$. I.e. it
suffices to show:

\[ \begin{array}{rcll}  tr(ag_i) & = & tr(g_ia) & \ \ \ \ a \in {\cal H}_{n+1}, \ i=1,\ldots, n 
\\ [.1in]
                      tr(at) & = & tr(ta)  & \ \ \ \ a \in {\cal H}_{n+1}. \end{array} \]

\noindent By Theorem 4, $a$ is of one of the following types:  

\vspace{.1in}
\, \, i) $a=w_{n-1}$

\vspace{.1in}
\,  ii) $a=w_{n-1}g_ng_{n-1} \ldots g_i$

\vspace{.1in}
iii) $a=w_{n-1}g_ng_{n-1} \ldots g_{i+1}{t'_i}^k, k\in \ZZ$ respectively $\ZZ_d$

\vspace{.1in}
iv) $a=w_{n-1}{t'_n}^k, k\in \ZZ$ respectively  $\ZZ_d$, where $w_{n-1}$ is some word in  ${\cal
H}_n$. 

\bigbreak 
\noindent If $a=w_{n-1}$ and $x=t$ or  $x=g_i$ for $ i=1,\ldots, n-1 $ the assumption
holds from the induction step, whilst for $x=g_n$ it follows from (3$'$) that $tr(w_{n-1}g_n)=
z \, tr(a) =  tr(g_nw_{n-1})$. 

\bigbreak
\noindent If $a$ is of type ii) or of type iii) and $x=t$ or  $x=g_i$ for $ i=1,\ldots, n-1 $ we
apply the same reasoning as above using rule  (3$'$). So we have  to check still the cases where 
$a=w_{n-1}g_ng_{n-1} \ldots g_i$ or $a=w_{n-1}g_ng_{n-1} \ldots g_{i+1}{t'_i}^k$ and $x=g_n$,
i.e.

\[ \begin{array}{rcl}
tr(w_{n-1}g_n \ldots g_ig_n) & = & tr(g_nw_{n-1}g_n \ldots g_i)  \\ [.1in]  
tr(w_{n-1}g_n \ldots g_{i+1}{t'_i}^kg_n) & = & tr(g_n w_{n-1}g_n \ldots g_{i+1}{t'_i}^k)
\end{array} \ \ \ \ \  (*) \]

\noindent If $a$ is of type iv) and $x=t$ or  $x=g_i$ for $ i=1,\ldots, n-1 $ we have to check:

\[ \begin{array}{rcl}
tr(w_{n-1}{t'_n}^k t) & = & tr(t w_{n-1}{t'_n}^k) \\ [.1in]  
tr(w_{n-1}{t'_n}^k g_i) & = & tr(g_i w_{n-1}{t'_n}^k) 
\end{array} \ \ \ \ \ \ \ (**) \]

\noindent Finally, if $a$ is of type iv)  and $x=g_n$  we have to check:

\[ tr(w_{n-1}{t'_n}^k g_n) \ = \ tr(g_n w_{n-1}{t'_n}^k) \ \ \ \ \ \ \ (***) \]

 Before checking $(*), (**)$ and $(***)$  we need the following: 
 
\begin{lem}
The function $tr$ satisfies the following stronger version of rule 4):

\[ (4') \ \ tr(x{t_n'}^ky) = s_k \, tr(xy),\] 

\noindent for any $x, y \in {\cal H}_n, \ k \in \ZZ \mbox{ respectively }\ZZ_d. $

\end{lem}
 
\begin{pf} It suffices to prove (4$'$) for the case that $y$ is of the form $y=y_1 t^\lambda
y_2$, where $y_1$is a product of the $g_i$'s for $ i=1,\ldots, n-1,  \ \lambda \in \ZZ 
\mbox{ respectively }\ZZ_d$ and $y_2$ an arbitrary word in ${\cal H}_n$. Indeed we have: 

\bigbreak
\noindent $tr(x{t_n'}^ky) = tr(x\underline{{t_n'}^k y_1} t^\lambda
y_2)\stackrel{Lemma 1,(vi)}{=} tr(xy_1{t_n'}^k t^\lambda y_2)$
    
\vspace{.1in} 
$=tr(xy_1g_n\ldots g_1 t^k{g_1}^{-1}\underline{{g_2}^{-1} \ldots {g_n}^{-1} t^\lambda}
y_2)\stackrel{Lemma 1,(vi)}{=}$
    
\vspace{.1in} 
$=tr(xy_1g_n\ldots \underline{ g_1 t^k{g_1}^{-1}t^\lambda }{g_2}^{-1} \ldots {g_n}^{-1})=A$

\vspace{.1in}
\noindent The latter underlined expression says that we have to consider four possibilities 
depending on $k, \lambda$ being positive or negative. We show here the case where both  
$k, \lambda$ are positive. The rest are proved completely analogously. For $k, \lambda$ 
positive, Lemma 5,(i) says:   

\[ \begin{array}{lcl}
g_1 t^k {g_1}^{-1} t^\lambda & = &  t^\lambda g_1 t^k {g_1}^{-1}  + 
(q^{-1} -1) \, [t^{\lambda -1} g_1 t^{k+1} + \cdots +  g_1 t^{k+\lambda}] \\ [.1in]  
& + &  (1-q^{-1}) \, [t^k g_1 t^\lambda + \cdots + t^{k+\lambda -1} g_1 t ]. 
\end{array}\]

\noindent We substitute then in $A$ to obtain:
 
\vspace{.1in} 
\noindent $ A= tr(xy_1\underline{g_n\ldots g_2 t^\lambda} g_1 t^k {g_1}^{-1} \ldots
{g_n}^{-1}y_2) $
 
\vspace{.1in}
$+ (q^{-1} -1) \, [tr(xy_1\underline{g_n\ldots g_2 t^{\lambda -1}} g_1
\underline{t^{k+1}{g_2}^{-1} \ldots {g_n}^{-1}}y_2) + \cdots  $
 
\vspace{.1in}
$+ tr(xy_1g_n\ldots g_1 \underline{t^{k+\lambda} {g_2}^{-1} \ldots {g_n}^{-1}}y_2)]  $   
 
\vspace{.1in}
$+ (1-q^{-1}) \, [tr(xy_1\underline{g_n\ldots g_2 t^k} g_1 \underline{t^\lambda {g_2}^{-1} \ldots
{g_n}^{-1}}y_2) + \cdots  $

\vspace{.1in}
$+ tr(xy_1\underline{g_n\ldots g_2 t^{k+\lambda -1}} g_1 \underline{t {g_2}^{-1} \ldots
{g_n}^{-1}}y_2)] \stackrel{Lemma\, 1,(vi)}{=}$
 
\vspace{.1in}
 $= tr(xy_1 t^\lambda {t'_n}^k y_2) $
 
\vspace{.1in}
$ + (q^{-1} -1) \, [tr(xy_1 t^{\lambda -1} \underline{ g_n\ldots g_1 {g_2}^{-1} \ldots {g_n}^{-1}}
t^{k+1} y_2) + \cdots  $
 
\vspace{.1in}
$+ tr(xy_1\underline{ g_n\ldots g_1 {g_2}^{-1} \ldots {g_n}^{-1}} t^{k+\lambda} y_2)] $   
 
\vspace{.1in}
$ +(1-q^{-1}) \, [tr(xy_1 t^k \underline{g_n\ldots g_1 {g_2}^{-1} \ldots {g_n}^{-1}}t^\lambda y_2)
+ \cdots  $

\vspace{.1in}
$+ tr(xy_1 t^{k+\lambda -1} \underline{g_n\ldots g_1 {g_2}^{-1} \ldots{g_n}^{-1}} t y_2)]
\stackrel{Lemma\, 1,(iii)}{=}$
 
\vspace{.1in}
 $= tr(xy_1 t^\lambda {t'_n}^k y_2)  $
 
\vspace{.1in}
$+ (q^{-1} -1) \, [tr(xy_1 t^{\lambda -1} {g_1}^{-1} \ldots {g_{n-1}}^{-1}  \underline{g_n}\ldots
g_1  t^{k+1} y_2) + \cdots $
 
\vspace{.1in}
$+ tr(xy_1{g_1}^{-1} \ldots {g_{n-1}}^{-1} \underline{g_n}\ldots g_1 t^{k+\lambda} y_2)]  $   
 
\vspace{.1in}
$+ (1-q^{-1}) \, [tr(xy_1 t^k {g_1}^{-1} \ldots  {g_{n-1}}^{-1} \underline{g_n}\ldots g_1
t^\lambda y_2) + \cdots $

\vspace{.1in}
$ +tr(xy_1 t^{k+\lambda -1} {g_1}^{-1} \ldots {g_{n-1}}^{-1} \underline{g_n}\ldots g_1 t y_2)]
\stackrel{(3')}{=}$

\vspace{.1in}
 $= tr(xy_1 t^\lambda {t'_n}^k y_2) + (q^{-1} -1) z \, [tr(xy_1 t^{\lambda +k} y_2) + (1-q^{-1})z
\, [tr(xy_1 t^{k+\lambda} y_2) $

\vspace{.1in}
$= tr(xy_1 t^\lambda {t'_n}^k y_2). $

\end{pf} 

\noindent The relations $(**)$ follow now immediately from Lemma 7, since:

\vspace{.1in}
\noindent $tr(w_{n-1}\underline{{t'_n}^k} g_i)  \stackrel{(4')}{=} s_k \, tr(w_{n-1} g_i)
\stackrel{induction \, step}{=}s_k \, tr(g_i w_{n-1}) = tr(g_i w_{n-1}{t'_n}^k), $  

\vspace{.1in}
\noindent  for all $i<n$, and similarly for $x=t$.  

\bigbreak
We next show $(*)$  for $a=w_{n-1}g_n \ldots g_i$. The case $a=w_{n-1}g_n \ldots g_{i+1}{t'_i}^k$
is shown similarly. On the one hand we have: 

\vspace{.1in}
 $tr(w_{n-1}g_ng_{n-1}\underline{ \ldots g_ig_n})= tr(w_{n-1}\underline{g_ng_{n-1}g_n}g_{n-2}
\ldots g_i)$

\vspace{.1in}
 $=tr(w_{n-1}g_{n-1}\underline{g_n}g_{n-1}g_{n-2}\ldots g_i)\stackrel{(3')}{=} 
z\, tr(w_{n-1}\underline{{g_{n-1}}^2}g_{n-2}\ldots g_i) $

\vspace{.1in}
 $= (q-1)z\, tr(w_{n-1}g_{n-1}\ldots g_i)+ qz\, tr(w_{n-1}g_{n-2}\ldots g_i). $

\vspace{.1in}
\noindent On the other hand in order to calculate $tr(g_nw_{n-1}g_n \ldots g_i)$ we examine the
different possibilities for $w_{n-1}$:

\vspace{.1in}
\noindent -- If $w_{n-1}\in {\cal H}_{n-1},$ then $tr(\underline{g_nw_{n-1}}g_n\ldots g_i) 
 =tr(w_{n-1}\underline{{g_n}^2}g_{n-1}\ldots g_i) $

\vspace{.1in}
$=(q-1)z\, tr(w_{n-1}g_{n-1}\ldots g_i) + qz\, tr(w_{n-1}g_{n-2}\ldots g_i).$

\vspace{.1in}
\noindent -- If $w_{n-1}=bg_{n-1}c$, where $b, \, c \in {\cal H}_{n-1}$, then
$tr(\underline{g_nb}g_{n-1}\underline{cg_n}g_{n-1}\ldots g_i) $

\vspace{.1in}
 $=tr(bg_{n-1}\underline{g_n}g_{n-1} c g_{n-1}\ldots g_i) \stackrel{(3')}{=}
z\, tr(b\underline{{g_{n-1}}^2} c g_{n-1}\ldots g_i)$

\vspace{.1in}
 $= (q-1)z\, tr(bg_{n-1} c g_{n-1}\ldots g_i)+ qz\, tr(b c g_{n-1}\ldots g_i)$

\vspace{.1in}
 $=(q-1)z\, tr(bg_{n-1} c g_{n-1}\ldots g_i)+ qz^2\, tr(b c g_{n-2}\ldots g_i)$

\vspace{.1in}
 $=(q-1)z\, tr(bg_{n-1} c g_{n-1}\ldots g_i)+ qz\, tr(b g_{n-1}c g_{n-2}\ldots g_i)$

\vspace{.1in}
 $=(q-1)z\, tr(w_{n-1} g_{n-1}\ldots g_i)+ qz\, tr(w_{n-1} g_{n-2}\ldots g_i).$

\vspace{.1in}
\noindent -- Finally, if $w_{n-1}=b{t'_{n-1}}^k$, where $b, \in {\cal H}_{n-1}$, then 

\vspace{.1in}
 $tr(\underline{g_nb}{t'_{n-1}}^kg_n \ldots g_i)=tr(bg_n{t'_{n-1}}^k\underline{g_n} \ldots g_i)$

\vspace{.1in}
 $=q\, tr(b \underline{{t'_n}^k g_{n-1}} \ldots g_i) + 
(q-1)\, tr(b \underline{ g_n}{t'_{n-1}}^kg_{n-1} \ldots g_i)\stackrel{(4'),(3')}{=}$

\vspace{.1in}
 $=qz\, tr(b{t'_{n-1}}^k g_{n-2}\ldots g_i) + (q-1)z\, tr(b {t'_{n-1}}^k g_{n-1} \ldots g_i)$

\vspace{.1in}
 $=qz\, tr(w_{n-1} g_{n-2}\ldots g_i) + (q-1)z\, tr(w_{n-1} g_{n-1} \ldots g_i).$

\begin{note} \rm The relations $(*)$ imply that  $tr(xg_nyg_n)=tr(g_nxg_ny)$ for any $x,\, y 
\in {\cal H}_n$. 

\end{note}
 
\noindent It remains now to show $(***)$. On the one hand we have:

\vspace{.1in}
 $tr(w_{n-1} \underline{{t'_n}^kg_n}) \stackrel{Lemma \, 1,(vi)}{=}
tr(w_{n-1} \underline{g_n}{t'_{n-1}}^k) \stackrel{(3')}{=}z\, tr(w_{n-1} {t'_{n-1}}^k).$

\vspace{.1in}
\noindent On the other hand in order to calculate $tr(g_nw_{n-1}{t'_n}^k)$ we examine the
different possibilities for $w_{n-1}$:

\vspace{.1in}
\noindent -- If $w_{n-1}\in {\cal H}_{n-1},$ then $tr(\underline{g_nw_{n-1}}{t'_n}^k) = 
tr(w_{n-1}\underline{{g_n}^2}{t'_{n-1}}^k{g_n}^{-1})$

\vspace{.1in}
$ =(q-1)\, tr(w_{n-1}\underline{{t'_n}^k}) + q\, tr(w_{n-1}{t'_{n-1}}^k\underline{{g_n}^{-1}}) =
(q-1)\, tr(w_{n-1}{t'_{n-1}}^k) $

\vspace{.1in}
$+z\,  tr(w_{n-1}{t'_{n-1}}^k) + (1-q)\,  tr(w_{n-1}{t'_{n-1}}^k)=z\,  tr(w_{n-1}{t'_{n-1}}^k).$

\vspace{.1in}
\noindent -- If $w_{n-1}=ag_{n-1}b,$ where $a,\, b \in {\cal H}_{n-1},$ then 

\vspace{.1in}
$tr(\underline{g_na}g_{n-1}b\underline{{t'_n}^k}) =  
tr(a\underline{g_ng_{n-1}g_n}b{t'_{n-1}}^k\underline{{g_n}^{-1}} $

\vspace{.1in}
$ =q^{-1}\, tr(\underline{ag_{n-1}}g_n \underline{g_{n-1}b{t'_{n-1}}^k}g_n) + 
(q^{-1}-1)\, tr(ag_{n-1}\underline{g_n} g_{n-1}b{t'_{n-1}}^k) = $

\vspace{.1in}
(applying  Note 2  for $x=ag_{n-1}$ and $y=g_{n-1}b{t'_{n-1}}^k)$

\vspace{.1in}
$= q^{-1}\, tr(\underline{g_na}g_{n-1}g_n g_{n-1}b{t'_{n-1}}^k) + 
(q^{-1}-1)z\, tr(a{g_{n-1}}^2 b{t'_{n-1}}^k)$

\vspace{.1in}
$ = q^{-1}\, tr(ag_{n-1}\underline{g_n}{g_{n-1}}^2 b{t'_{n-1}}^k) + 
(q^{-1}-1)z\, tr(a{g_{n-1}}^2 b{t'_{n-1}}^k) $

\vspace{.1in}
$= q^{-1}z(q^2-q+1)\, tr(\underline{ag_{n-1}b}{t'_{n-1}}^k) + 
q^{-1}zq(q-1)\, tr(a b{t'_{n-1}}^k)$

\vspace{.1in}
$ + (q^{-1}-1)z(q-1)\, tr(a\underline{g_{n-1}}b{t'_{n-1}}^k) +
(q^{-1}-1)zq\, tr(ab{t'_{n-1}}^k) =z\, tr(w_{n-1}{t'_{n-1}}^k). $

\bigbreak
\noindent Before proving the last case we need to deform the expression ${t'_{n-1}}^l{t'_n}^k$.
Indeed we have:

\vspace{.1in}
\noindent ${t'_{n-1}}^l{t'_n}^k = g_{n-1}\ldots g_1t^l\underline{{g_1}^{-1}\ldots
{g_{n-1}}^{-1}g_n\ldots  g_1}t^k{g_1}^{-1}\ldots {g_n}^{-1}$

\vspace{.1in}
 $=g_{n-1}\ldots g_1t^l\underline{g_n\ldots g_2}g_1\underline{{g_2}^{-1}\ldots{g_n}^{-1}}         
t^k{g_1}^{-1}\ldots {g_n}^{-1}$

\vspace{.1in}
 $=(g_{n-1}g_n)\ldots (g_1g_2)t^lg_1 t^k ({g_2}^{-1}{g_1}^{-1})\ldots ({g_n}^{-1}{g_{n-1}}^{-1})
\underline{{g_n}^{-1}}$

\vspace{.1in}
 $=(g_{n-1}g_n)\ldots (g_1g_2)t^lg_1 t^k {g_1}^{-1}   ({g_2}^{-1}{g_1}^{-1})\ldots
({g_n}^{-1}{g_{n-1}}^{-1})$

\vspace{.1in}
 $=q^{-1}\, (g_{n-1}g_n)\ldots (g_1g_2) \underline{t^lg_1 t^k g_1}  
({g_2}^{-1}{g_1}^{-1})\ldots ({g_n}^{-1}{g_{n-1}}^{-1})$

\vspace{.1in}
$+(q^{-1}-1)\, (g_{n-1}g_n)\ldots (g_1g_2) t^lg_1 t^k g_1  
({g_2}^{-1}{g_1}^{-1})\ldots ({g_n}^{-1}{g_{n-1}}^{-1})\stackrel{Lemma \, 5,(i)}{=}$

\vspace{.1in}
 $=q^{-1}\, (g_{n-1}g_n)\ldots (g_1g_2) \underline{g_1}t^kg_1 t^l   
({g_2}^{-1}{g_1}^{-1})\ldots ({g_n}^{-1}{g_{n-1}}^{-1})$

\vspace{.1in}
$+(1-q^{-1})\, [(g_{n-1}g_n)\ldots (g_1g_2) tg_1 t^{k+l-1}  
({g_2}^{-1}{g_1}^{-1})\ldots ({g_n}^{-1}{g_{n-1}}^{-1})+\cdots $

\vspace{.1in}
$+(g_{n-1}g_n)\ldots (g_1g_2) t^lg_1 t^k  
({g_2}^{-1}{g_1}^{-1})\ldots ({g_n}^{-1}{g_{n-1}}^{-1})] $

\vspace{.1in}
$+(q^{-1}-1)\, [(g_{n-1}g_n)\ldots (g_1g_2) t^kg_1 t^l  
({g_2}^{-1}{g_1}^{-1})\ldots ({g_n}^{-1}{g_{n-1}}^{-1})+\cdots $

\vspace{.1in}
$+(g_{n-1}g_n)\ldots (g_1g_2) t^{k+l-1}g_1 t  
({g_2}^{-1}{g_1}^{-1})\ldots ({g_n}^{-1}{g_{n-1}}^{-1})] $

\vspace{.1in}
$+(q^{-1}-1)\, (g_{n-1}g_n)\ldots (g_1g_2) t^lg_1 t^k  
({g_2}^{-1}{g_1}^{-1})\ldots ({g_n}^{-1}{g_{n-1}}^{-1}) $

\vspace{.1in}
 $=q^{-1}\, g_ng_{n-1}\ldots g_1 t^k \underline{g_n\ldots g_1{g_2}^{-1}\ldots {g_n}^{-1}}
t^l  {g_1}^{-1}\ldots {g_{n-1}}^{-1}$

\vspace{.1in}
$+(1-q^{-1})\, [g_{n-1}\ldots g_1 t\underline{g_n\ldots g_1{g_2}^{-1}\ldots {g_n}^{-1}} t^{k+l-1}  
{g_1}^{-1}\ldots {g_{n-1}}^{-1} +\cdots $

\vspace{.1in}
$+g_{n-1}\ldots g_1 t^l\underline{g_n\ldots g_1{g_2}^{-1}\ldots {g_n}^{-1}} t^k  
{g_1}^{-1}\ldots {g_{n-1}}^{-1}] $

\vspace{.1in}
$+(q^{-1}-1)\, [g_{n-1}\ldots g_1t^k\underline{g_n\ldots g_1{g_2}^{-1}\ldots {g_n}^{-1}} t^l  
{g_1}^{-1}\ldots {g_{n-1}}^{-1} +\cdots $

\vspace{.1in}
$+g_{n-1}\ldots g_1 t^{k+l-1}\underline{g_n\ldots g_1{g_2}^{-1}\ldots {g_n}^{-1}} t  
{g_1}^{-1}\ldots {g_{n-1}}^{-1}] $

\vspace{.1in}
$+(q^{-1}-1)\, g_{n-1}\ldots g_1 t^l\underline{g_n\ldots g_1{g_2}^{-1}\ldots {g_n}^{-1}} t^k  
{g_1}^{-1}\ldots {g_{n-1}}^{-1} $

\vspace{.1in}
 $=q^{-1}\, g_n\ldots g_1 t^k {g_1}^{-1}\ldots {g_{n-1}}^{-1}l\underline{g_n}\ldots g_1 t^l  
{g_1}^{-1}\ldots{g_{n-1}}^{-1}$

\vspace{.1in}
$+(1-q^{-1})\, [g_{n-1}\ldots g_1 t{g_1}^{-1}\ldots {g_{n-1}}^{-1}g_n\ldots g_1 t^{k+l-1}  
{g_1}^{-1}\ldots {g_{n-1}}^{-1} +\cdots $

\vspace{.1in}
$+g_{n-1}\ldots g_1 t^{l-1}{g_1}^{-1}\ldots {g_{n-1}}^{-1}g_n\ldots g_1 t^{k+1}  
{g_1}^{-1}\ldots {g_{n-1}}^{-1}] $

\vspace{.1in}
$+(q^{-1}-1)\, [g_{n-1}\ldots g_1t^k{g_1}^{-1}\ldots {g_{n-1}}^{-1}g_n\ldots g_1 t^l  
{g_1}^{-1}\ldots {g_{n-1}}^{-1} +\cdots $

\vspace{.1in}
$+g_{n-1}\ldots g_1 t^{k+l-1}{g_1}^{-1}\ldots {g_{n-1}}^{-1}g_n\ldots g_1 t  
{g_1}^{-1}\ldots {g_{n-1}}^{-1}] $

\vspace{.1in}
 $=g_n\ldots g_1 t^k {g_1}^{-1}\ldots g_n^{-1}g_{n-1}\ldots g_1 t^l  
{g_1}^{-1}\ldots{g_{n-1}}^{-1}$

\vspace{.1in}
 $+(1-q^{-1})\, g_n\ldots g_1 t^{k+l} {g_1}^{-1}\ldots {g_{n-1}}^{-1}$

\vspace{.1in}
$+(1-q^{-1})\, [g_{n-1}\ldots g_1 t{g_1}^{-1}\ldots {g_{n-1}}^{-1}g_n\ldots g_1 t^{k+l-1}  
{g_1}^{-1}\ldots {g_{n-1}}^{-1} +\cdots $

\vspace{.1in}
$+g_{n-1}\ldots g_1 t^{l-1}{g_1}^{-1}\ldots {g_{n-1}}^{-1}g_n\ldots g_1 t^{k+1}  
{g_1}^{-1}\ldots {g_{n-1}}^{-1}] $

\vspace{.1in}
$+(q^{-1}-1)\, [g_{n-1}\ldots g_1t^k{g_1}^{-1}\ldots {g_{n-1}}^{-1}g_n\ldots g_1 t^l  
{g_1}^{-1}\ldots {g_{n-1}}^{-1} +\cdots $

\vspace{.1in}
$+g_{n-1}\ldots g_1 t^{k+l-1}{g_1}^{-1}\ldots {g_{n-1}}^{-1}g_n\ldots g_1 t  
{g_1}^{-1}\ldots {g_{n-1}}^{-1}].$

\vspace{.1in}
\noindent Notice that with applying the other cases of Lemma 5 we obtain analogous results.

\bigbreak
\noindent -- If, finally, $w_{n-1}=b{t'_{n-1}}^l,$ where $b \in {\cal H}_{n-1},$ we have: 
$tr(\underline{g_n}b \underline{{t'_{n-1}}^l{t'_n}^k})  $

\vspace{.1in}
 $=tr(b\underline{{g_n}^2}g_{n-1}\ldots g_1 t^k {g_1}^{-1}\ldots g_n^{-1}g_{n-1}\ldots g_1 t^l  
{g_1}^{-1}\ldots{g_{n-1}}^{-1})$

\vspace{.1in}
 $+(1-q^{-1})\, tr(b\underline{{g_n}^2}g_{n-1}\ldots g_1 t^{k+l}{g_1}^{-1}\ldots {g_{n-1}}^{-1})$

\vspace{.1in}
$+(1-q^{-1})\, [tr(bg_ng_{n-1}\ldots g_1 t{g_1}^{-1}\ldots {g_{n-1}}^{-1}\underline{g_n}\ldots
g_1 t^{k+l-1}{g_1}^{-1}\ldots {g_{n-1}}^{-1})  $

\vspace{.1in}
$+\cdots+ tr(bg_ng_{n-1}\ldots g_1 t^{l-1}{g_1}^{-1}\ldots {g_{n-1}}^{-1}\underline{g_n}\ldots g_1
t^{k+1}   {g_1}^{-1}\ldots {g_{n-1}}^{-1})] $

\vspace{.1in}
$+(q^{-1}-1)\, [tr(bg_ng_{n-1}\ldots g_1t^k{g_1}^{-1}\ldots {g_{n-1}}^{-1}\underline{g_n}\ldots g_1 t^l  
{g_1}^{-1}\ldots {g_{n-1}}^{-1}) +\cdots $

\vspace{.1in}
$+tr(g_{n-1}\ldots g_1 t^{k+l-1}{g_1}^{-1}\ldots {g_{n-1}}^{-1}\underline{g_n}\ldots g_1 t  
{g_1}^{-1}\ldots {g_{n-1}}^{-1})] $

\vspace{.1in}
 $=(q-1)\, tr(b{t'_n}^k{t'_{n-1}}^l) + q\, tr(b{t'_{n-1}}^k \underline{{g_n}^{-1}}{t'_{n-1}}^l)+
(1-q^{-1})[(q-1)z+q]\, tr(b{t'_{n-1}}^{k+l})$

\vspace{.1in}
$+(1-q^{-1})\, [q\, tr(bt'_n{t'_{n-1}}^{k+l-1}) +(q-1)z\, tr(b{t'_{n-1}}^{k+l}) +\cdots $

\vspace{.1in}
$+q\, tr(b{t'_n}^{l-1}{t'_{n-1}}^{k+1})+(q-1)z\, tr(b{t'_{n-1}}^{k+l})] $

\vspace{.1in}
$+(q^{-1}-1)\, [q\, tr(b{t'_n}^k{t'_{n-1}}^l) +(q-1)z\, tr(b{t'_{n-1}}^{k+l}) +\cdots $

\vspace{.1in}
$+q\, tr(b{t'_n}^{k+l-1}t'_{n-1})+(q-1)z\, tr(b{t'_{n-1}}^{k+l})] $

\vspace{.1in}
 $=(q-1)s_ks_l\, tr(b) + q[q^{-1}z+(q^{-1}-1)]\, tr(b{t'_{n-1}}^{k+l})+$

\vspace{.1in}
$[(q+q^{-1}-2)z+(q-1)]\, tr(b{t'_{n-1}}^{k+l}) + (q^{-1}-1)(q-1)z\, tr({t'_{n-1}}^{k+l})tr(b)$

\vspace{.1in}
$+(q-1)s_1\, tr({t'_{n-1}}^{k+l-1})tr(b) +\cdots +(q-1)s_{l-1}\, tr({t'_{n-1}}^{k+1})tr(b) $

\vspace{.1in}
$+(1-q)\, tr({t'_n}^k)s_l\, tr(b) +\cdots +(1-q)\, tr({t'_n}^{k+l-1})s_1\, tr(b).$

\vspace{.1in}
 And since $tr({t'_n}^i)=tr({t'_{n-1}}^i)$ in all algebras ${\cal H}_n$, we conclude that  

\vspace{.1in}
$tr(g_nb {t'_{n-1}}^l{t'_n}^k) = z\, tr(b {t'_{n-1}}^{k+l}) = z\, tr(b {t'_{n-1}}^l{t'_{n-1}}^k)
= z\, tr(w_{n-1}{t'_{n-1}}^k).$

\bigbreak
\noindent The proof of Theorem 6 is now concluded.
\end{pf}

As already mentioned in the Introduction, we can define $tr$ with so few rules, because 
 the elements ${t}^k,\ldots, {t'_i}^k$ in rule 4) are all conjugate, and  this reflects the fact 
that  $B_{1,n}$ splits as a semi-direct product of the classical braid group $B_n$ and of its
free  subgroup $ P_{1,n} $ generated precisely by the elements $t, t'_1, \ldots, t'_{n-1}$: 
 \ $ B_{1,n} = P_{1,n} \semi \, B_n.$ 

\bigbreak
Note that if $k\in \ZZ_2$ we are in the case of the classical Iwahori-Hecke algebras of type  ${\cal
B}$, and from the above construction we recover the trace given in  \cite{L1,L2}. Moreover, if a word 
$x \in {\cal H}_n$ does not contain any $t$'s (that is, if $x$ is an element of the 
Iwahori-Hecke algebra of type  ${\cal A}_n$), then $tr(x)$ can be computed using only rules 1),
2), and 3) of Theorem 6, and in this case $tr$ agrees with Ocneanu's trace (cf. \cite{J}).

\begin{rem} \rm 
A word seen as an element of different ${\cal B}$-type Hecke algebras will aquire in
principle different values for the different traces. This difference consists in substituting
-- if necessary -- the parameters $s_i$ according to the defining relation $(\spadesuit)$ of 
${\cal H}_n(q,d): \ t^d=a_{d-1}t^{d-1}+\cdots + a_0$.  So, in  ${\cal H}_n(q,d)$ we have: 
$tr({t'_n}^k)=s_k$ for $k\in \ZZ_d$ and $tr({t'_n}^d)=a_{d-1}s_{d-1}+\cdots + a_0.$

\vspace{.1in}
For example in ${\cal H}_n(q,\infty)$ and
in ${\cal H}_n(q,d)$ for $d>5$ we have $tr(t^5)=s_5$. 

\vspace{.1in}
In ${\cal H}_n(q,5)$ is 
$tr(t^5)=a_4s_4+\cdots + a_0,$ whilst in ${\cal H}_n(q,3)$ is

\vspace{.1in} 
$tr(t^5)=({a_2}^3+ 2a_1a_2+ a_0)s_2 + ({a_1}^2+ a_1{a_2}^2+ a_0a_2)s_1 + 
(a_0a_1+ a_0{a_2}^2).$

\end{rem}

\noindent In order to calculate the trace of a word in ${\cal H}_n$ we bring it to the canonical
form of Theorem 5 applying at the same time the rules of the trace. As an example  we calculate
below $tr(g_2g_1t^3g_1^{-1}g_3g_2g_3)$. We have:   

\vspace{.1in}
\noindent
$tr(g_2g_1t^3g_1^{-1}\underline{g_3g_2g_3})=tr(g_2g_1t^3g_1^{-1}g_2\underline{g_3}g_2)=z \,
tr(g_2g_1t^3g_1^{-1}\underline{{g_2}^2})$

\vspace{.1in}
 $=z(q-1) \, tr(g_2g_1t^3g_1^{-1}\underline{g_2}) + zq \, tr(g_2g_1t^3g_1^{-1})$

\vspace{.1in}
 $=z(q-1)q \, tr({t'_2}^3) + z(q-1)^2 \, tr(\underline{g_2}g_1t^3g_1^{-1}) + zq \,
tr(\underline{g_2}g_1t^3g_1^{-1})$

\vspace{.1in}
 $=q(q-1)z \, tr({t'_2}^3) + z^2(q^2-q+1) \, tr({t'_1}^3). $

\section{Invariants of knots in the solid torus}

The aim of this section is to construct all analogues of the 2-variable Jones polynomial
homfly-pt) for oriented knots in the solid torus derived from the cyclotomic and
generalized Hecke algebras of type ${\cal B}$,  using their Markov equivalence and the Markov
traces constructed in Theorem 6. All knots/links will be assumed to be oriented, and we shall say
`knots' for both knots and links. 
\bigbreak

As mentioned in the Introduction the elements of the braid groups $B_{1,n}$, which we call
`mixed braids', are represented geometrically by braids in $n+1$ strands in $S^3$, which keep the
first strand fixed. The closure of a mixed braid  represents a knot inside the oriented solid
torus, $ST$, where the fixed strand represents the complementary solid torus in $S^3$, and the
next $n$ numbered strands represent the  knot in $ST$. Below we illustrate a mixed braid in 
$B_{1,5}$ and a knot in $ST$. 

$$\vbox{\picill3.3inby1.2in(canon5)  }$$

\noindent Moreover, it has been well-understood that all  knots in $ST$ may be
represented by mixed braids, and their isotopy in $ST$ is reflected by equivalence classes of
braids in  $\bigcup_{n=1}^{\infty}B_{1,n}$ through the following: 

\begin{th}{\rm (cf.\ \cite{L2}, Theorem~3.)}\\
Let $L_1$,  $L_2$  be  two oriented  links  in $ST$   and  ${\beta}_1$,
${\beta}_2$  be  mixed   braids  in  $\bigcup_{n=1}^{\infty}  B_{1,n}$
corresponding to $L_1$, $L_2$. Then $L_1$ is  isotopic to $L_2$ in $ST$
if and   only  if    ${\beta}_1$  is  equivalent  to   ${\beta}_2$  in
$\bigcup_{n=1}^{\infty} B_{1,n}$  under  equivalence generated by  the
braid relations together with the following two moves:
\begin{itemize}
\item[(i)] Conjugation: If $\alpha,\beta \in B_{1,n}$ then $\alpha \sim \beta ^{-1}\alpha\beta$.
\item[(ii)] Markov moves: If $\alpha \in B_{1,n}$ then $\alpha \sim \alpha{\sigma_n}^{\pm 1}\in
B_{1,n+1}$.  \end{itemize}
\end{th}

Let  now $\pi$  denote    the canonical  quotient  map \  $B_{1,n}
\longrightarrow {\cal H}_n$ \ given in Definition~1, and consider the trace constructed in Theorem
6 for a specified algebra ${\cal H}_n$.   Then   a braid  in $B_{1,n}$  can  be mapped 
through  $tr \circ \pi$ to  an expression in  the variables $q, u_1^{\pm 1}, \ldots, u_d^{\pm 1},
\ldots, z, (s_k), \ k\in \ZZ$ respectively $\ZZ_d$. Let also $\widehat{\alpha}$ 
denote the knot obtained by closing the mixed braid $\alpha$. 
Theorem 7 combined with Theorem 6 say that in order to obtain a knot invariant ${\cal X}$ in $ST$
from  any specified trace of Theorem 6 we  have to normalize first $g_i$ to $\sqrt{\lambda} g_i$
so that 

\[ tr(a (\sqrt{\lambda} g_n))=tr(a({(\sqrt{\lambda} g_n)}^{-1})) \ \ \mbox{for \ } a \in 
{\cal H}_n.\]

\noindent This normalization has been done in \cite{L2}, (5.1), where   Jones's
normalization of Ocneanu's trace (cf.  \cite{J}) was followed, and it yields 

\[ \lambda:=\frac{z+1-q}{qz}, \  \ z:=\frac{1-q}{q\lambda-1}.\]

\noindent Then we have to normalize $tr$ so that  

\[{\cal X} (\widehat{\alpha})={\cal X}(\widehat{\alpha{\sigma_n}}) = {\cal X}(\widehat{\alpha
\sigma_n^{-1}}).\] 

\noindent Let finally $A$ be  the field  of rational functions  over
$\QQ$  in   indeterminates $\sqrt{\lambda}, \sqrt{q},$ $a_{d-1}, \ldots, a_0, \ldots, (s_k),
\ k\in \ZZ$ respectively $\ZZ_d$. (The reason for having square root of $q$ becomes clear in
the recursive formula $\dagger$ below.)
 Then the normalizations result the following

\begin{defn}{\rm
(cf.\  \cite{L2}, Definition 1.)}  \rm  For $\alpha$, $tr$, $\pi$ as above    let  

\[{\cal X}_{\widehat{\alpha}}  =    {\cal X}_{\widehat{\alpha}}(q, a_{d-1}, \ldots, a_0,
\sqrt{\lambda},  s_1, s_2,   \ldots)  :=  
\Bigl[-\frac{1-\lambda q}{\sqrt{\lambda}(1-q)} \Bigr]^{n-1}
(\sqrt{\lambda})^e\, tr(\pi(\alpha)), \] 

\noindent where $e$ is the exponent sum of the $\sigma_i$'s that appear in  $\alpha$.  (Note that  
the $t'_i$'s do  not    affect   the   estimation   of   $e$, so they can be ignored.) Then
${\cal X}_{\hat{\alpha}}$   depends    only  on     the   isotopy  class   of the mixed  knot 
$\widehat{\alpha}$, which represents an oriented  knot in $ST$. (For example, in ${\cal H}_n(q,d)$
and for $k\in \ZZ_d$ we have: $\alpha=t^k$, then ${\cal X}_{\hat{\alpha}}=s_k$.) 
\end{defn}

Note that if a knot in $ST$ can be enclosed in a 3-ball then it may be seen as a knot in  $S^3$ 
and there exists a mixed braid representative, $\alpha$, which does not contain $t'_i$'s. Then 
${\cal X}_{\hat{\alpha}}$ has  the   same value as the  2-variable Jones polynomial 
(homfly-pt) as given in \cite{J}, Definition~6.1.  On the lower level of ${\cal H}_n(q,Q) \ {\cal X}$
 yields the invariants constructed  in \cite{L2}, Section 5 and \cite{GL}, Section 5.

\begin{rem} \rm 
Note furthermore that one could also define ${\cal H}_n(q,d)$ as a quotient of $B_{1,n}$ by
sending the generator $t$ of $B_{1,n}$ to $t^{-1}$ of ${\cal H}_n(q,d)$. Then the traces and
the knot invariants in $ST$ constructed above exhaust the whole range of such constructions
related to {\it all possible} Hecke and Hecke-related algebras of type ${\cal B}$.  \end{rem}

\noindent {\bf On recursive formulae:} We shall now show how to interpret the above in terms of knot
diagrams, and how to calculate  alternatively  the above knot invariants in $ST$ by applying
recursive {\it skein relations} and {\it initial conditions} on the mixed link diagrams.
 Let $L_+ ,   L_-  , L_0$ be   oriented  mixed link diagrams that   are
identical, except in one crossing, where they are as depicted below:

$$\vbox{\picill3inby0.75in(canon6)  }$$

\noindent With analogous reasoning as in \cite{L2}, (5.2) (cf. also \cite{J}) the defining
quadratic relation of ${\cal H}_n$ induces the invariant ${\cal X}$ to satisfy the following
recursive linear formula, which is  the
well-known skein rule used   for   the evaluation of   the   homfly-pt
polynomial.

 \[\frac{1}{\sqrt{q}\sqrt{\lambda}}\, {\cal X}_{L_+} - \sqrt{q}\sqrt{\lambda}\,
 {\cal X}_{L_-}  = (\sqrt{q}-\frac{1}{\sqrt{q}})\, {\cal X}_{L_0} \ \ \ \ \ \ \ \dagger  \]

\noindent In the case of  ${\cal H}_n(q, \infty)$ there is no other skein relation that ${\cal
X}$ satisfies. 

\noindent In the case of  ${\cal H}_n(q, d)$, let  $M_{d}, M_{d-1}, \ldots, M_0$  be oriented
mixed link diagrams   that are identical, except in the regions depicted below:

$$\vbox{\picill3inby1in(canon7)  }$$

\noindent Using conjugation we may assume that  $M_{d}=\widehat{\alpha\, {t'_i}^d}, 
M_{d-1}=\widehat{\alpha\, {t'_i}^{d-1}}, \ldots, M_{0}=\widehat{\alpha}$ for some 
$\alpha\in B_{1,n}$. And so by Lemma 1, (viii) we obtain: 

\[tr(\pi(\alpha\, {t'_i}^d))  =  a_{d-1}\, tr(\pi(\alpha\, {t'_i}^{d-1}))   +\cdots + 
a_0\, tr(\pi(\alpha)),  \]

\noindent If we multiply now the above equation by  

\[ \Bigl[ - \frac{1-\lambda q}{\sqrt{\lambda}(1- q)} \Bigr]^{n-1}\, {(\sqrt{\lambda})}^e \] 

\noindent we obtain the following skein relation for ${\cal X}$:

 \[{\cal X}_{\widehat{\alpha {t'_i}^d}} =  a_{d-1}\, {\cal X}_{\widehat{\alpha {t'_i}^{d-1}}}
+ \cdots + a_0\,{\cal X}_{M_0} \ \ \ \ \ \ \ \ddagger  \]
 
\noindent (compare with Remark 3). We  next find the
initial conditions that are also needed for  evaluating ${\cal X}$ for any knot diagram  in $ST$
using the skein relations $\dagger$ and $\ddagger$. Clearly 

\[{\cal X}_{unknot}=1\]

\noindent should be one of them. Recall now the canonical basis of ${\cal H}_{n+1}$ given in
Theorem 5. With appropriate changes of crossings (using the quadratic relations for the $g_i$'s)
this basis yields a canonical enumeration of descending diagrams related to $B_{1,n+1}$. Let now 
$\alpha$ be such a 
descending diagram. Applying $tr$ on $\alpha$ means geometrically that we close the braid
$\alpha$ and we apply the Markov moves. Using Rule (4) we extract and re-insert 
$tr({t'_i}^k)$ so as to obtain:

\[{\cal X}_{\widehat{\alpha}} =  
{\cal X}_{\widehat{{t'_{i_1}}^{k_1}{t'_{i_2}}^{k_2}\ldots {t'_{i_r}}^{k_r}} }.  \]
 
\noindent This provides the second set of initial conditions, namely the values of ${\cal X}$ at
all links consisting of stucks of loops of different twists with same orientation around the
`axis' solid torus. If ${\cal X}$ is derived by the cyclotomic Hecke algebra ${\cal H}_n(q,d)$ 
the number of twists of each loop cannot exceed $d-1$. In the case of ${\cal H}_n(q,\infty)$ the
number of twists is arbitrary. We illustrate below an example of a descending diagram with the
starting point at the top of the last strand, the basic link $t^4t_1{t_2}^{-1}$ and the projection
of $t^3$ on a punctured disc.  

$$\vbox{\picill4.5inby1.4in(canon8)  }$$

We conclude with some remarks. 

\bigbreak 
\noindent {\bf Remarks}  (i) \ On the level of ${\cal H}_n(q,\infty), {\cal X}$ is defined by
all initial conditions (with unrestricted number of twists) and only by the first skein rule.
Therefore the set of mixed links of the form  $\widehat{{t'_{i_1}}^{k_1}{t'_{i_2}}^{k_2}\ldots
{t'_{i_r}}^{k_r}},$ for $k_1, \ldots, k_r \in \ZZ$ forms the  basis of the 3rd skein module of
the solid torus. Thus 
 the result of J.Hoste and M.Kidwell in \cite{HK}, and of V.Turaev in \cite{T} is recovered
with this method. If ${\cal X}$ is derived by the cyclotomic  Hecke algebras ${\cal H}_n(q,d)$ the
set of mixed links of the form  $\widehat{{t'_{i_1}}^{k_1}{t'_{i_2}}^{k_2}\ldots
{t'_{i_r}}^{k_r}},$ for $k_1, \ldots, k_r \in \ZZ_d$ forms the  basis of the corresponding
submodule of the 3rd skein module of $ST$. In \cite{P} the algebra ${\cal H}_n(q,\infty)$
has been studied independently and the corresponding $ST$-invariant has been constructed using
similar methods.

\bigbreak 
\noindent (ii) \ If on the level of ${\cal H}_n(q,\infty)$ we use the skein rule 
 \[\frac{1}{t}\, {\cal Y}_{L_+} - t\, {\cal Y}_{L_-}  = 
(\sqrt{t}-\frac{1}{\sqrt{t}})\, {\cal Y}_{L_0} \] 

\noindent instead of $\dagger$, and the initial conditions  ${\cal Y}_{unknot}=1$ and 
${\cal Y}_{\widehat{t}} = s$ we obtain an analogue of the Jones polynomial for oriented knots in
the {\it oriented} $ST$. If $ST$ is {\it unoriented} we have to allow an extra isotopy move for
knots in $ST$, namely to flip over the diagram around the $x$-axis, where the knot diagram is
projected on a punctured disc.  The invariant $\cal Y$ is preserved under the flipping over move, so
$\cal Y$ is the analogue of the Jones polynomial in the {\it orientable} $ST$. For details and for the
Kauffman bracket approach of this invariant  see  \cite{HP}.  

\bigbreak 
\noindent (iii) \ The invariant ${\cal X}$ related to  ${\cal H}_n(q,\infty)$ is the
appropriate one for extending the results to the lens spaces $L(p,1)$. The combinatorial setup
is similar to the one for $ST$, only the Markov braid equivalence  includes one more move,
which reflects the surgery description of $L(p,1)$. So, in order to construct a homfly-pt
analogue for knots in  $L(p,1)$ or, equivalently, in order to compute for  $L(p,1)$  the 3rd skein
module and its quotients  we have to normalize the $ST$-invariants further so that 

\[{\cal X}_{\widehat{\alpha}} =  {\cal X}_{sl(\widehat{\alpha})}, \ \mbox{for \ }\alpha \in
B_{1,n}, \]

\noindent  for all possible slidings of $\alpha$. This is the subject of [S.~Lambropoulou,
J.~Przytycki, Hecke algebra approach to the skein module of lens spaces, in preparation].

\bigbreak 
\noindent (iv) \ Analogous combinatorial setup, Markov braid equivalence and 
braid structures in arbitrary c.c.o. 3-manifolds and knot complements has already been  done in
\cite{LR},[S.~Lambropoulou, Braid structures in 3-manifolds, to appear in JKTR]. Therefore it is
possible in principle to extend such algebraic constructions to other 3-manifolds, by means of
constructing appropriate quotient algebras and Markov traces on them, followed by appropriate
normalizing, in order to derive knot invariants.

 
\bigbreak
\noindent {\sc S.L.: 3--5 Bunsenstrasse, Mathematisches Institut, G\"{o}ttingen Universit\"{a}t,
  37073 G\"{o}ttingen, Germany. E-mail: sofia@cfgauss.uni-math.gwdg.de}

\end{document}